\documentclass[11pt]{amsart}

\usepackage{amsmath}
\usepackage{amssymb}
\usepackage{mathrsfs}
\usepackage{color}

\newtheorem{theorem}{Theorem}[section]
\newtheorem{claim}[theorem]{Claim}

\newtheorem{lemma}[theorem]{Lemma}

\newtheorem{proposition}[theorem]{Proposition}
\newtheorem{corollary}[theorem]{Corollary}

\theoremstyle{definition}
\newtheorem{definition}[theorem]{Definition}

\newtheorem{question}[theorem]{Question}

\theoremstyle{remark}
\newtheorem{remark}[theorem]{Remark}

\newcount\skewfactor
\def\mathunderaccent#1#2 {\let\theaccent#1\skewfactor#2
\mathpalette\putaccentunder}
\def\putaccentunder#1#2{\oalign{$#1#2$\crcr\hidewidth
\vbox to.2ex{\hbox{$#1\skew\skewfactor\theaccent{}$}\vss}\hidewidth}}

\def\smallbox#1{\leavevmode\thinspace\hbox{\vrule\vtop{\vbox
   {\hrule\kern1pt\hbox{\vphantom{\tt/}\thinspace{\tt#1}\thinspace}}
   \kern1pt\hrule}\vrule}\thinspace}

\def\l{{\langle}}
\def\r{{\rangle}}

\newcommand{\lusim}[1]{\smash{\underset{\raisebox{1.2pt}[0cm][0cm]{$\sim$}}
{{#1}}}}

\newcommand{\cf}{{\rm cf}}
\def\s{\subseteq}
\newcommand{\Col}{\mathop{\mathrm{Col}}}
\DeclareMathOperator{\id}{id}

\DeclareMathOperator{\crit}{crit}

\title[Galvin's property at large cardinals and partition calculus]{Galvin's property at large cardinals and an application to partition calculus}

\author{Tom Benhamou}
\address{Departement of Mathematics, Statistics, and Computer Science,  University of Illinois at Chicago, USA}
\email{tomb@uic.edu}

\author{Shimon Garti}
\address{Einstein Institute of Mathematics,
 The Hebrew University of Jerusalem,
 Jerusalem 91904, Israel}
\email{shimon.garty@mail.huji.ac.il}

\author{Alejandro Poveda}
\address{Center of Mathematical Sciences and Applications, Harvard University, Cambridge MA, 02138, USA}
\email{alejandro@cmsa.fas.harvard.edu}

\thanks{}

\begin{document}
\let\labeloriginal\label
\let\reforiginal\ref
\def\ref#1{\reforiginal{#1}}
\def\label#1{\labeloriginal{#1}}

\begin{abstract}
In the first part of this paper we explore the possibility for a very large cardinal $\kappa$ to carry a $\kappa$-complete ultrafilter without Galvin's property. In this context we prove the consistency of every ground model $\kappa$-complete ultrafilter extends to a non-Galvin one. Oppositely, it is also consistent that every ground model $\kappa$-complete ultrafilter extends to a $P$-point ultrafilter, hence to another one satisfying Galvin's property. 
Finally, we apply this property to obtain consistently new instances of the classical problem in partition calculus  $\lambda\rightarrow(\lambda,\omega+1)^2$. 
\end{abstract}

\maketitle

\section{Introduction}
Let $\mathscr{F}$ be a filter over a regular uncountable cardinal $\kappa$.
We say that \emph{Galvin's property} holds for $\mathscr{F}$ (in symbols, ${\rm Gal}(\mathscr{F},\kappa,\kappa^+)$) if every family $\langle C_\gamma\mid \gamma<\kappa^+\rangle \subseteq\mathscr{F}$\footnote{Note that, in the absence of \textsf{AC}, it is important to emphasis that a family is a \textit{sequence} in which we allow repetitions. } admits a subfamily $\langle C_{\gamma_i}\mid i<\kappa\rangle$ with the property that $\bigcap\{C_{\gamma_i}\mid i<\kappa\}\in\mathscr{F}$.
In the 1970's, Galvin proved that if $\kappa=\kappa^{<\kappa}>\aleph_0$ then ${\rm Gal}(\mathscr{F},\kappa,\kappa^+)$ is true whenever $\mathscr{F}$ is normal.
The statement and the proof were  published in a paper by Baumgartner, Hajnal and Mate \cite{MR0369081}. 

The motivation of this paper came from an open problem which appeared in \cite{bgs}.
In that work it is shown that, consistently, there is a $\kappa$-complete ultrafilter over a measurable cardinal $\kappa$ which fails to satisfy the Galvin property.
One should keep in mind the fact that if $\kappa$ is measurable then every normal filter $\mathscr{F}$ satisfies the Galvin property $\mathrm{Gal}(\mathscr{F},\kappa,\kappa^+)$.
Thus, the main result of \cite{bgs} shows that $\kappa$-completeness differs from normality in terms of implying Galvin's property.
On the other hand, it is consistent that $\kappa$ is measurable and every $\kappa$-complete ultrafilter is Galvin.
This can be demonstrated in Solovay's inner model $L[\mathscr{U}]$, as shown in \cite{MR4393795}.
Howe\-ver, inner models are limited with their tolerance to large cardinals.
It was asked in \cite{bgs} whether it is consistent for a supercompact cardinal $\kappa$ that 
every $\kappa$-complete ultrafilter $\mathscr{U}$ over $\kappa$ satisfies $\mathrm{Gal}(\mathscr{U},\kappa,\kappa^+)$. 
In the first part of this paper we investigate the possibility of very large cardinals carrying $\kappa$-complete filters (ultrafilters) that fail to satisfy Galvin's property. 
In \S2.1 we exhibit  a generic extension where $\kappa$ is supercompact and every $\kappa$-complete ground model ultrafilter $\mathscr{U}$ over $\kappa$ extends to an ultrafilter $\mathscr{U}^*$ for which  $\mathrm{Gal}(\mathscr{U}^*,\kappa,\kappa^+)$ fails (see Proposition~\ref{ExtendingGalvinWithSupercompacts}). Shortly after, we show that this construction is amenable to preserve even stronger large cardinals, such as $C^{(n)}$-extendibles and Vop\v{e}nka's Principle (Proposition~\ref{Largercadinals}). Continuing in this vein, we present a result of an opposite nature. 
Namely, in Theorem~\ref{MakingEverythingNormal} we construct a generic extension where $\kappa$ is supercompact and every $\kappa$-complete ground model ultrafilter $\mathscr{U}$ extends to a $\kappa$-complete ultrafilter $\mathscr{U}^*$ that is \emph{Rudin-Keisler} equivalent to a normal one. In particular, all of these ultrafilters $\mathscr{U}^*$ do satisfy Galvin's property (Proposition~\ref{IdkappathenGalvin}).  
The reader may have noticed that this  is perhaps a too harsh way to convert an arbitrary $\kappa$-complete ultrafilter into a Galvin one. During the rest of the section we present alternative strategies to achieve the same configuration without such dramatic changes. Our first attempt takes place in Theorem~\ref{ExtendFilters} where we employ iterations of \emph{Generalized Mathias forcing} to show that every $\kappa$-complete filter $\mathscr{U}$ extends to a $\kappa$-complete filter $\mathscr{U}^*$ for which $\mathrm{Gal}(\mathscr{U}^*,\kappa,\kappa^+)$ holds. Section 2.1 is then culminated with our main result, which builds upon previous work of Gitik and Shelah \cite{MR1632081}. More specifically, in  Theorem~\ref{GitikShelahTheorem} we replace the previous iteration by a more sophisticated one also involving Generalized Mathias forcing. This iteration was devised by Gitik and Shelah and here it is adapted to our current purposes. As an outcome we obtain the consistency of a supercompact cardinal $\kappa$ with every ground model $\kappa$-complete ultrafilter $\mathscr{U}$ extending to a $\kappa$-complete ultrafilter $\mathscr{U}^*$ which is a $P$-point. In particular, $\mathrm{Gal}(\mathscr{U}^*,\kappa,\kappa^+)$ holds (see Proposition~\ref{P point is Galvin}).

\smallskip

\smallskip

The study of Galvin's property led us to the area of partition calculus. 
In \S3 of the paper we address a classical question about ordinary partition relations. An exquisite theorem of Shelah establishes that if $\lambda>\cf(\lambda)=\kappa>\aleph_0$ and $2^\kappa<\lambda$ then $\lambda\rightarrow(\lambda,\omega+1)^2$  \cite{MR2494318}.  We prove that the same partition relation follows upon replacing the cardinal arithmetic assumption by an appropriate instance of Galvin's property. This is true in general, but it is particularly interesting  under \textsf{AD}. Indeed, under this assumption Galvin's property holds for a wide class of cardinals (i.e. \emph{Boldface $\mathsf{GCH}$ cardinals}). In fact, if $\kappa$ is measurable and $\kappa=\cf(\lambda)<\lambda$ is a limit of measurable cardinals then $\lambda\rightarrow(\lambda,\omega+1)^2$.
This gives an answer to \cite[Question 11.4]{MR795592} in the context of \textsf{AD}. 
For details, see Theorem~\ref{thm881ad} and the subsequent discussion. 

Finally, it  must be said that  our Galvin-like assumption is trivial when $2^\kappa<\lambda$, and forceable when $2^\kappa>\lambda$. Thus, we get more positive instances of $\lambda\rightarrow(\lambda,\omega+1)^2$. We believe, however, that the  relation $\lambda\nrightarrow(\lambda,\omega+1)^2$ is consistent as well. Actually, our result pinpoints which instances of Galvin's property should be violated in order to force this negative partition relation.

\smallskip

Our notation is mostly standard.
If $\kappa=\cf(\kappa)>\aleph_0$ then $\mathscr{D}_\kappa$ denotes the club filter over $\kappa$.
If $\kappa=\cf(\kappa)<\lambda$ then $S^\lambda_\kappa=\{\delta\in\lambda\mid \cf(\delta)=\kappa\}$. The arrow symbol $\lambda\rightarrow (\alpha,\beta)^2$ is a shorthand for the following statement:  for every $f:[\lambda]^{2}\rightarrow 2$ either there is a $0$-monochoromatic subsets of $\lambda$ of order type $\alpha$ or a $1$-monochromtic subset of $\lambda$ of order-type $\beta$.
We say that $\binom{\alpha}{\beta}\rightarrow\binom{\gamma}{\delta}$ iff for every $c:\alpha\times\beta\rightarrow\{0,1\}$ there are $A\in[\alpha]^\gamma,B\in[\beta]^\delta$ for which $c\upharpoonright(A\times{B})$ is constant.
We use $\Theta$ to denote $$\sup\{\alpha\mid \text{There exists a mapping from ${}^\omega\omega$ \emph{onto} $\alpha$}\}.$$
We employ the Jerusalem forcing notation, thus $p\leq{q}$ means that $q$ is stronger than $p$.
For background in partition calculus we refer the reader to \cite{MR795592} and \cite{MR3075383}.


\section{Galvin's property at very large cardinals with and without choice} 
\subsection{Galvin's property at large cardinals}
In \cite{MR4393795} the following result is proved: It is consistent that $\kappa$ is a measurable cardinal and every $\kappa$-complete  ultrafilter $\mathscr{U}$ over $\kappa$ satisfies $\mathrm{Gal}(\mathscr{U},\kappa,\kappa^+)$.
The proof strategy is based on analyzing  Solovay's inner model $L[\mathscr{U}]$, where a complete classification of the $\sigma$-complete ultrafilters over $\kappa$ is available. The key observation is that in this inner model every $\sigma$-complete ultrafilter over $\kappa$ is Rudin-Keisler equivalent to a finite power of the normal measure $\mathscr{U}$. Since these ultrafilters do satisfy Galvin's property one concludes that $\mathrm{Gal}(\mathscr{V},\kappa,\kappa^+)$ holds for every $\kappa$-complete ultrafilter $\mathscr{V}\in L[\mathscr{U}]$.   This phenomenon suggests the following  question:  How about those (very) large cardinals for which there is no  available canonical inner model?  The epitome of this is supercompactness.

By work of the first two authors together with S. Shelah \cite{bgs} it is consistent that a supercompact cardinal $\kappa$ carries a  $\kappa$-complete ultrafilter $\mathscr{U}$ which extends the club filter and $\neg \mathrm{Gal}(\mathscr{U},\kappa,\kappa^+)$. Shortly after, the first author together with M. Gitik \cite{OnPrikryandCohen} improved this result by showing that just a measurable cardinal suffices to obtain such an ultrafilter $\mathscr{U}$.

The forthcoming proposition is a spin-off of the above-mentioned result in the context of general $\kappa$-complete ultrafilters: 
\begin{proposition}\label{ExtendingGalvinWithSupercompacts}
 Assume that the $\mathsf{GCH}$ holds and that $\kappa$ is a measurable cardinal. Then the following is true in the generic extension of \cite[Theorem 2.6]{OnPrikryandCohen}: 
 Every $\kappa$-complete (not necessarily normal) ultrafilter $\mathscr{U}$ of the ground model extends  to a $\kappa$-complete ultrafilter $\mathscr{U}^*$ such that $\neg \mathrm{Gal}(\mathscr{U}^*,\kappa,\kappa^+)$. 
 
 In addition, if $\kappa$ was supercompact then it remains so in the  extension.
\end{proposition}
\begin{proof}
The sought model is the generic extension by the Easton support iteration $\l\mathbb{P}_\alpha,\lusim{\mathbb{Q}}_\beta\mid \alpha\leq\kappa+1, \beta\leq\kappa\r$ such that for $\alpha\leq\kappa$, $\lusim{\mathbb{Q}}_\alpha$ is trivial unless $\alpha$ is inaccessible, in which case it is a $\mathbb{P}_\alpha$-name for $\mathrm{Add}(\alpha,\alpha^+)$. This iteration preserves supercompactness (see e.g. \cite[Theorem~11.1]{MR2768691}). 

Let ${\mathscr{U}}\in V$ be a  $\kappa$-complete ultrafilter. Let us verify that we can adjust the argument in \cite{OnPrikryandCohen} to encompass non-normal ultrafilters. We will follow the notation from the original proof,  considering the elementarity embeddings
$$j_1:=j_{{\mathscr{U}}}:V\rightarrow M_{\mathscr{U}}=:M_1, \ j_2:=j_{\mathscr{U}^2}:V\rightarrow M_{\mathscr{U}^2}=:M_2$$
 $$ k:M_1\rightarrow M_{2}, \ j_{2}=k\circ j_1$$
where $k$ is simply the ultrapower embedding defined in $M_\mathscr{U}$ using the ultrafilter $j_{1}(\mathscr{U})$. Let
 $G:=G_\kappa*g_\kappa$ be $V$-generic for $\mathbb{P}_\kappa*\lusim{\mathbb{Q}}_\kappa$. The argument that these embeddings can be lifted in $V[G]$ does not require normality and remains unaltered. Thus, we form  $j_1\subseteq j_1^*:V[G]\rightarrow M_1[j^*_1(G)]$, $k\subseteq k^*:M_1[j_1^*(G)]\rightarrow M_2[j_2^*(G)]$ and $j_2\subseteq j_2^*:=k^*\circ j_1^*$ such that:
 \begin{enumerate}
     \item for every $\alpha\in j_1``\kappa^+$, $f_{\kappa_2,k(\alpha)}(\kappa_1)=1$.
     \item for every $\alpha\in \kappa_1\setminus j_1``\kappa^+, f_{\kappa_2,k(\alpha)}(\kappa_1)=0.$
     \item $f_{\kappa_2,\kappa_1}(\kappa_1)=\kappa$.
 \end{enumerate}
 Since we are just dealing with non-normal ultrafilters we need to alter the values of the generic $f_2$ at $\delta^*:=[\id]_{j_1(\mathscr{U})}$,  the generator of the second ultrapower. Also, we need to eliminate the generator of the first ultrapower $\delta:=[\id]_\mathscr{U}$:
 \begin{enumerate}
     \item for every $\alpha\in j_1``\kappa^+$, $f_{\kappa_2,k(\alpha)}(\delta^*)=1$.
     \item for every $\alpha\in \kappa_1\setminus j_1``\kappa^+$, $f_{\kappa_2,k(\alpha)}(\delta^*)=0$
     \item $f_{\kappa_2,\delta^*}(\delta^*)=\delta$.
 \end{enumerate}
Notice that the amount of  coordinates that were altered is small. In particular,
 the counting/genericity arguments of \cite[Lemma 2.7]{OnPrikryandCohen} relying on \textsf{ZFC} still go through. Next, derive in $V[G]$ the  ultrafilter generated by  $j^*_1$ and  $[\id]_\mathscr{U}$,
$$\mathscr{U}^*:=\{X\subseteq \kappa\mid [\id]_\mathscr{U}\in j^*_1(X)\}$$
Note that $\mathscr{U}\subseteq \mathscr{U}^*$. Finally, let 
\begin{equation*}
	\mathscr{W}:=\{X\subseteq\kappa\mid [\id]_{j_1(\mathscr{U})}\in j^*_2(X)\}\in V[G].\qedhere
\end{equation*}

Let us prove that $\mathscr{W}$ witnesses the statement of the theorem:
\begin{claim}\label{final claim}
$\mathscr{W}$ is a $\kappa$-complete ultrafilter over $\kappa$ such that:
\begin{enumerate}
    \item  $\mathscr{U}\subseteq \mathscr{W}$.
\item $\neg \mathrm{Gal}(\mathscr{W},\kappa,\kappa^+)$.
\end{enumerate}
\end{claim}
\begin{proof}[Proof of claim]
 (1): If $A\in \mathscr{U}$ then $j_1(A)\in j_1(\mathscr{U})$, hence $[\id]_{j_1(\mathscr{U})}\in j_2(A)$ and thus $A\in \mathscr{W}$.

 (2): Let us define the witness. For each $\alpha<\kappa^+$ let
 $$A_\alpha:=\{\nu<\kappa\mid f_{\kappa,\alpha}(\nu)=1\}$$
 then $$j^*_2(A_\alpha)=\{\beta<\kappa_2\mid f_{\kappa_2,j_2(\alpha)}(\beta)=1\}.$$ Since $j_2(\alpha)=k(j_1(\alpha))$, our modifications of the generic give $$f_{\kappa_2,j_2(\alpha)}([\id]_{j_1(\mathscr{U})})=1,$$ hence $[\id]_{j_1(\mathscr{U})}\in j^*_2(A_\alpha)$. Finally  $A_\alpha\in \mathscr{W}$ by definition of $\mathscr{W}$.
 Before proving the  failure of the Galvin property, let us denote by $j_\mathscr{W}:V[G]\rightarrow M_\mathscr{W}$ the ultrapower embedding by $\mathscr{W}$ and $k_\mathscr{W}:M_\mathscr{W}\rightarrow M^*_2$ defined by $k_\mathscr{W}([f]_\mathscr{W}):=j^*_2(f)([id]_{j_1(\mathscr{U})})$ the factor map satisfying $k_\mathscr{W}\circ j_\mathscr{W}=j^*_2$. 
 
 We show that $k_\mathscr{W}$ is onto, hence the identity, and thus $j^*_2=j_\mathscr{W}$. In effect, if $A\in M_2[j^*_2(G)]$ then there is a name $\lusim{A}\in M_2$ with $A=(\lusim{A})_{j_2^*(G)}$. Since $j_2$ is the second ultrapower by $\mathscr{U}$, there is $f:[\kappa]^2\rightarrow V$ such that $j_2(f)([\id]_\mathscr{U},[\id]_{j_1(\mathscr{U})})=\lusim{A}$. By \L\"{o}s theorem, we can assume that $f(\alpha,\beta)$ is a $\mathbb{P}_{\kappa+1}$-name for every $(\alpha,\beta)\in [\kappa]^2$. In $V[G]$ let $f^*(\alpha)=(f(f_{\kappa,\alpha}(\alpha),\alpha))_G$. 
 
 Then, $$k_\mathscr{W}([f^*]_W)=j_2^*(f^*)([\id]_{j_1(\mathscr{U})})=(j_2(f)(f_{\kappa_2,[\id]_{j_1(\mathscr{U})}}([\id]_{j_1(\mathscr{U})}),[\id]_{j_1(\mathscr{U})}))_{j_2^*(G)}$$
 $$=j_2(f)([\id]_\mathscr{U},[\id]_{j_1(\mathscr{U})})_{j_2^*(G)}=(\lusim{A})_{j_2^*(G)}=A$$
 Let
 $\langle A_{\alpha_i}\mid i<\kappa\rangle$ be any subfamily of length $\kappa$ and $\kappa\leq \eta<[\id]_\mathscr{W}=[\id]_{j_1(\mathscr{U})}$.
 Denote $j_\mathscr{W}(\l A_{\alpha_i}\mid i<\kappa\r):=\l A'_{\alpha'_i}\mid i<j_\mathscr{W}(\kappa)\r.$

 Pick any  $\kappa\leq\eta< [\id]_{\mathscr{U}}<j_2(\kappa)$, then $\eta\notin j_1``\kappa^+$ and also $\alpha''_\eta\notin j_1``\kappa^+$, where $\alpha''_\eta$ is the first image of the $\{\alpha_i\mid i<\kappa\}$.  Moreover $k(\alpha''_\eta)=\alpha'_{k(\eta)}$ and by definition, $f_{\kappa_2,\alpha'_{\eta}}([\id]_{j_1(\mathscr{U})})=0$ and $[\id]_{j_1(\mathscr{U})}\notin A'_{\eta}$. Hence
  $$[\id]_{j_1(\mathscr{U})}\notin \bigcap \{A'_{\alpha'_i}\mid i<\kappa_2\}=j_2^*(\bigcap_{i<\kappa}A_{\alpha_i})$$
  Hence $\bigcap_{i<\kappa}A_{\alpha_i}\notin \mathscr{W}$.
\end{proof}

\end{proof}

Continuing with our original discussion one may ask if the conclusion of Proposition~\ref{ExtendingGalvinWithSupercompacts} is compatible with large cardinals stronger than supercompactness. As argued in \cite{Bag, BagPov}, the natural model-theoretic strengthe\-ning of supercompactness is $C^{(n)}$-extendibility. Fix $n<\omega$. A cardinal $\kappa$ is called \emph{$C^{(n)}$-extendible} if for every $\lambda>\kappa$ there is $\theta\in\mathrm{Ord}$ and an elementary embedding $j\colon V_\lambda\rightarrow V_\theta$ with $\crit(j)=\kappa$,  $j(\kappa)>\lambda$ and $V_{j(\kappa)}\prec_{\Sigma_n} V$.\footnote{Recall that $V_\eta\prec_{\Sigma_n} V$ is a shorthand for the following statement: for every $\bar{a}\in V_\eta^{<\omega}$ and every $\Sigma_n$ formula $\varphi(\bar{x})$ in the language of set theory, $V_\eta\models \varphi(\bar{a})$ iff $V\models \varphi(\bar{a})$.} 
 The classical notion of extendibility (see \cite[\S23]{Kan}) coincides with $C^{(n)}$-extendibility whenever $n=1$. However, when $n\geq 2$ the first $C^{(n)}$-extendible is far above, and has stronger large-cardinal-properties, than the first extendible cardinal. In addition, $C^{(n)}$-extendibility do entail a proper hierarchy of cardinals \cite{Bag}. The culmination of this hierarchy is 
 the category-theoretic axiom known as \emph{Vop\v{e}nka's Principle} (\textsf{VP}) \cite[p. 335]{Kan}. In effect, it was shown by Bagaria that \textsf{VP} is equivalent to the existence of a (proper class of) $C^{(n)}$-extendible, for all $n\geq 1$. We refer the reader to  \cite{Bag} for further details. 

Let us come back to the argument of Proposition~\ref{ExtendingGalvinWithSupercompacts}. If $\kappa$ is an extendible cardinal performing our iteration $\mathbb{P}_{\kappa+1}$ 
will ruin extendibility of $\kappa$.\footnote{Actually, adding a single Cohen subset to $\kappa$ does it.}  Nevertheless, if one forces with $\mathrm{Add}(\alpha,\alpha^+)$ at every inaccessible cardinal  the situation changes completely. In \cite{BagPov} the authors develop a general theory of preservation of extendible cardinals under class-forcing iterations. Specifically, in \cite[\S8]{BagPov} it is shown that many classical class-forcing iterations (e.g., Jensen's iteration to force the \textsf{GCH}) do preserve extendible cardinals, as well as 
$C^{(n)}$-extendible cardinals and Vop\v{e}nka's Principle ($\mathrm{VP}$). 
\smallskip

The following proposition is an easy corollary of \cite[Theorem~8.4]{BagPov}:

\begin{proposition}\label{Largercadinals}
Assume that the $\mathsf{GCH}$ holds and that $\kappa$ is a $C^{(n)}$-extendible cardinal for some $n\geq 1$. Let $\mathbb{P}$ denote  the Easton support class iteration forcing with $\mathrm{Add}(\alpha,\alpha^+)$ at each inaccessible cardinal.

Then, the following hold in $V^{ \mathbb{P}}$:
\begin{enumerate}
    \item $\kappa$ is $C^{(n)}$-extendible;
    \item for every measurable cardinal $\lambda$ every $\lambda$-complete ultrafilter $\mathscr{U}\in V$ extends to a $\lambda$-complete ultrafilter $\mathscr{U}^*$ such that $\neg \mathrm{Gal}(\mathscr{U}^*,\lambda,\lambda^+).$
\end{enumerate}

In addition, if one assumes $\mathrm{VP}$ this is preserved in $V^{\mathbb{P}}$.
\end{proposition}
\begin{proof}
 Clause~(1) is an immediate consequence of  \cite[Theorem~8.4]{BagPov}. 
 For Clause~(2) we argue as follows. If  $\lambda_0$ stands for the first $V$-inaccessible cardinal then $\mathbb{P}$ admits a gap at $\lambda_0^{++}$.\footnote{I.e., $\mathbb{P}\simeq \mathbb{P}_1\ast \lusim{\mathbb{P}_2}$  where $|\mathbb{P}_1|<\lambda_0^{++}$ and $\Vdash_{\mathbb{P}_1}\text{$``\lusim{\mathbb{P}}_2$ is $\lambda^{++}_0$-distributive''}$.} Thus $\mathbb{P}$ does not create new measurable cardinals  in $V^{\mathbb{P}}$ \cite[Corollary~2]{HamkinsGap}. Let $\lambda$ be a $V$-measurable cardinal and $\mathscr{U}$ a $\lambda$-complete ultrafilter in the ground model. By Proposition~\ref{ExtendingGalvinWithSupercompacts}, $\mathbb{P}_{\lambda+1}$ forces that there is $\mathscr{U}^*\supseteq \mathscr{U}$ such that $\mathrm{Gal}(\mathscr{U}^*,\lambda,\lambda^+)$ fails. Clearly   $\mathbb{P}/\mathbb{P}_{\lambda+1}$ is forced to be $\lambda^{++}$-directed closed (actually more), hence it preserves that  $\mathscr{U}^*$ is a $\lambda$-complete ultrafilter over $\lambda$ witnessing $\neg \mathrm{Gal}(\mathscr{U}^*,\lambda,\lambda^+)$.
\end{proof}

Let us now show how to obtain the consistency of the opposite statement ``Every $\kappa$-complete ultrafilter $\mathscr{U}$ over a supercompact cardinal $\kappa$ extends to a $\kappa$-complete ultrafilter $\mathscr{U}^*$ satisfying $\mathrm{Gal}(\mathscr{U}^*,\kappa,\kappa^+)$''.  To this aim we will show how to turn a non (necessarily) Galvin ultrafilter $\mathscr{U}$ into another $\mathscr{U}^*$ that is \emph{Rudin-Keisler} equivalent to a normal one. By virtue of Galvin's theorem \cite{MR0369081} this ensures that $\mathscr{U}^*$ itself does satisfy $\mathrm{Gal}(\mathscr{U}^*,\kappa,\kappa^{+})$.


To show this we need a couple of preliminary observations. First, if $\mathscr{U}$ and $\mathscr{W}$ are  $\kappa$-complete ultrafilters over $\kappa$ and $\mathscr{U}\equiv_{\mathrm{RK}}\mathscr{W}$ then $\mathrm{Gal}(\mathscr{U},\kappa,\kappa^+)$ holds if and only if $\mathrm{Gal}(\mathscr{W},\kappa,\kappa^+)$ holds. Second:  
\begin{proposition}\label{IdkappathenGalvin}
If $\mathscr{U}$ is a $\kappa$-complete ultrafilter over $\kappa$ with $|[\id]_\mathscr{U}|=\kappa$ then $\mathscr{U}$ is Rudin-Keisler equivalent to a normal $\kappa$-complete ultrafilter. 

In particular, under the above conditions, $\mathrm{Gal}(\mathscr{U},\kappa,\kappa^+)$ holds.
\end{proposition}
\begin{proof}
Let  $\mathscr{U}_0$ denote the normal measure generated from $j:=j_\mathscr{U}$ and $\kappa$. 
 
 For each $\lambda<\kappa^+$ there is $f_\lambda\colon \kappa\rightarrow \kappa$ such that $j(f_\lambda)(\kappa)=\lambda$. We prove this by induction on $\lambda$. Suppose that  
 $\l f_\alpha\mid \alpha<\lambda<\kappa^+\r$ are defined and let $\l \lambda_i\mid i<\cf(\lambda)\r$  be cofinal in $\lambda$. Define $f_\lambda\colon \kappa\rightarrow \kappa$ as follows: $$f_\lambda(\alpha):=\sup_{i<\alpha}f_{\lambda_i}(\alpha).$$
 Note that $f_\lambda\colon \kappa\rightarrow\kappa$ due to the regularity of $\kappa$. Next, put $$j(\langle f_\beta\mid \beta<\lambda\rangle):=\langle f'_\beta\mid \beta<j(\lambda)\rangle, \ j(\l \lambda_i\mid i<\cf(\lambda)\r):=\l \lambda'_i\mid i<j(\cf(\lambda))\r.$$ Observe that $f'_{j(\alpha)}=j(f_\alpha)$ and $\lambda'_{j(\alpha)}=j(\lambda_\alpha)$. In particular,   $f'_\alpha=j(f_\alpha)$ and $\lambda'_\alpha=j(\lambda_\alpha)$ for every $\alpha<\kappa$. Hence, $$j(f_\lambda)(\kappa)=\sup_{i<\kappa}f'_{\lambda'_i}(\kappa)=\sup_{i<\kappa}j(f)_{j(\lambda_i)}(\kappa)=$$
 $$=\sup_{i<\kappa}j(f_{\lambda_i})(\kappa)=\sup_{i<\kappa}\lambda_i=\lambda$$
Thus, there is $f\colon \kappa\rightarrow \kappa$ such that $j(f)(\kappa)=[\id]_\mathscr{U}$, so that $\mathscr{U}\leq_{\mathrm{RK}} \mathscr{U}_0$. Also, it is well-known that normal filters are $\leq_{\mathrm{RK}}$-minimal (see~e.g. \cite[Proposition 2.6]{TomTreePrikry}), hence $\mathscr{U}\equiv_{\mathrm{RK}}\mathscr{U}_0$. For the in particular clause use our comments prior to the statement of the proposition.
\end{proof}

\begin{theorem}\label{MakingEverythingNormal}
Assume that the $\mathsf{GCH}$ holds and that $\kappa$ is a huge cardinal. 

Then, there is an inaccessible cardinal $\mu>\kappa$ and a generic extension of $V_\mu$ 
 where the following hold:
\begin{enumerate}
	\item $\kappa$ is supercompact;
	\item Every $\kappa$-complete {ultrafilter} $\mathscr{U}\in V$  extends to a $\kappa$-complete ultrafilter $\mathscr{U}^*$ that is Rudin-Keisler equivalent to a normal ultrafilter. In particular, $\mathrm{Gal}(\mathscr{U}^*,\kappa,\kappa^+)$ holds.
\end{enumerate}
\end{theorem}
\begin{proof}
	Let $j\colon V\rightarrow M$ be an elementary embedding witnessing that $\kappa$ is huge; namely, $\crit(j)=\kappa$ and $M^{j(\kappa)}\s M$.  Fix $\langle \mathscr{U}_\alpha\mid \alpha<2^{2^\kappa}\rangle$  an injective enumeration of the $\kappa$-complete ultrafilters over $\kappa$. For each $\alpha<2^{2^\kappa}$ note that  $j``\mathscr{U}_\alpha\in M$ and that  $j``\mathscr{U}_\alpha\in[j(\mathscr{U}_\alpha)]^{<j(\kappa)}$ hence, by $j(\kappa)$-completeness of $j(\mathscr{U}_\alpha)$ in $M$, we can find $\epsilon_\alpha\in \bigcap j``\mathscr{U}_\alpha.$ Clearly, $\epsilon_\alpha<j(\kappa)$ and 
	$$\mathscr{U}_\alpha\s\mathscr{U}^*_{\alpha,0}:= \{X\s \kappa\mid \epsilon_\alpha\in j(X)\}.$$
	Let $\lambda$ and $\mu$ be, respectively, the first inaccessible cardinals in the intervals $(\sup_{\alpha<2^{2^\kappa}}\epsilon_\alpha,  j(\kappa))$ and $(\lambda, j(\kappa))$.\footnote{This choice is possible as $j(\kappa)$ is a limit of inaccessibles.} Next, let $i\colon V\rightarrow N$ be the $\mu$-supercompact embedding derived from $j$; that is, the ultrapower embedding that arises from the measure $\{X\s\mathcal{P}_\kappa(\mu)\mid j``\mu\in j(X)\}$. Let $k\colon N\rightarrow M$ be the factor embedding between $j$ and $i$. Usual arguments show that $\crit(k)>\eta$, hence $\mathscr{U}^*_{\alpha,0}=\{X\s \kappa\mid \epsilon_\alpha\in i(X)\}$ for each $\alpha<2^{2^\kappa}$.
	
	Now, force over $V$ with \emph{Woodin's fast function forcing} $\mathbb{F}_\kappa$. By virtue of Lemma~1.10 in \cite{HamkinsLottery} we have that $i$ lifts to a $\mu$-supercompact embedding  $i\colon V[f]\rightarrow M[i(f)]$ such that $i(f)(\kappa)=\mu$. Notice that using the fast function $f\colon \kappa\rightarrow\kappa$ we can easily represent $\lambda$ as well: let $f^*\colon \kappa\rightarrow\kappa$ be defined as $\alpha\mapsto\sup\{\beta<f(\alpha)\mid \text{$\beta$ is inaccessible}\}$ and note that $i(f^*)(\kappa)=\lambda$.\footnote{Here we are implicitly assuming that $N[i(f)]$ is $\mu$-closed, hence it thinks that $\mu$ is inaccessible and that $\lambda$ is the first inaccessible below it.}
	
	Next, over $V[f]$, force with the two-step iteration $\mathbb{C}:=\mathbb{C}_\kappa\ast \lusim{\Col}(\kappa,<\lambda)$ where $\mathbb{C}_\kappa$ is the Easton-supported iteration defined as follows: for $\alpha<\kappa$, the $\alpha^{\mathrm{th}}$-stage of the iteration is trivial unless  $\alpha$ is inaccessible, $f``\alpha\s \alpha$ and $\alpha< f^*(\alpha)$, in which case it forces with $\lusim{\Col}(\alpha, {<}f^*(\alpha))$.
	\begin{claim}
		After forcing with $\mathbb{C}$ the embedding $i\colon V[f]\rightarrow N[f]$ lifts to a $\mu$-supercompact embedding $i^*\colon V[f\ast C]\rightarrow N[i(f\ast C)]$ in $V[f\ast C]$.
	\end{claim}
	\begin{proof}[Proof of claim]
	Denote $\bar{V}:=V[f]$ and $\bar{N}:=N[i(f)]$. Let $C:=C_\kappa\ast c\s \mathbb{C}$ generic over $V$. We can lift the embedding after forcing with $\mathbb{C}_\kappa$ to another $i\colon \bar{V}[C_\kappa]\rightarrow \bar{N}[C_\kappa\ast c\ast h]\s \bar{V}[C]$. There are two points here: first, the $\kappa^{\mathrm{th}}$-stage of the iteration from the perspective of $\bar{N}$ is $\Col(\kappa,{<}\lambda)^{\bar{V}}$;\footnote{Because $j(f)``\kappa\s \kappa$ and $j(f^*)(\kappa)=\lambda>\kappa$.} second, the tail forcing $i(\mathbb{C}_\kappa)/\mathbb{C}$ is trivial in the interval $(\kappa, \mu)$ because $i(f)(\kappa)=\mu$ and so the next closure point of $i(f)$ past $\kappa$ is $\geq (\mu^+)^{\bar{V}}$. 
	
	Finally, one can  lift $i$ after forcing with $\Col(\kappa,{<}\lambda)_{\bar{V}[C_\kappa]}$ to another embedding $i\colon \bar{V}[C]\rightarrow \bar{N}[C_\kappa\ast c\ast h \ast h']$. For this one uses the fact that $i``c\s \Col(i(\kappa),{<}i(\lambda))_{\bar{N}[C_\kappa\ast c\ast h]}$ is a directed set of conditions in $N[C_\kappa\ast c\ast h]$,  $|i``c|<i(\kappa)$ and $\Col(i(\kappa),{<}i(\lambda))_{\bar{N}[C_\kappa\ast c\ast h]}$ is $\mu$-directed-closed in the model  $N[C_\kappa\ast c\ast h]$. Standard arguments show that the resulting embedding witnesses $\mu$-supercompactness of $\kappa$.
	\end{proof}
	Working in  $V[f\ast C]$, for each $\alpha<(2^{2^\kappa})^V$ define $$\mathscr{U}^*_\alpha:=\{X\s\kappa\mid \epsilon_\alpha\in i^*(X)\}.$$
	Clearly, $\mathscr{U}^*_\alpha$ is a $\kappa$-complete ultrafilter satisfying $\mathscr{U}_\alpha\s \mathscr{U}^*_\alpha$. The point now is that $|[\id]_{\mathscr{U}^*_\alpha}|^{V[f\ast C]}\leq |\epsilon_\alpha|^{V[f\ast C]}=\kappa$ hence $\mathrm{Gal}(\mathscr{U}^*_\alpha,\kappa,\kappa^+)$ holds in $V[f\ast C]$. 
	
	Let $M:=V[f\ast C]_\mu$.  This is certainly a model of \textsf{ZFC} because $\mu$ remains inaccessible. Also, $M$ satisfies that $\kappa$ is supercompact. Finally, note that  $M=V_\mu[f\ast C]$. Since every $\kappa$-complete ultrafilter $\mathscr{U}$ over $\kappa$ from the ground model actually comes from $V_\mu$ we obtain Clause~(2) of the theorem. 
	\end{proof}
	In the model of Theorem~\ref{MakingEverythingNormal} our target cardinal $\kappa$ cannot be extendible. To make it so, one should perform a class-forcing iteration that is nice enough to carry the previous arguments. This suggests the following question: 
\begin{question}
		Is the statement of the previous theorem compatible with $\kappa$ being extendible or, more generally, $C^{(n)}$-extendible?
	\end{question}

In the proof of  Theorem~\ref{MakingEverythingNormal} we showed how to \emph{correct} ultrafilters that do not satisfy Galvin's property. The technique used for this purpose consisted of collapsing the generators of every $V$-ultrafilter to yet another generator of cardinality $\kappa$; namely, the normal generator. In that manner we accomplished our  \emph{correction} of non-Galvin ultrafilters by making them essentially \emph{minimal} (to wit, normal) from the Rudin-Keisler-perspective. This is, certainly,  a too harsh way to ensure Galvin's property in the final model.  

In the light of this one may wonder whether a similar Galvin-like configu\-ration is possible without \emph{trivializing} the relevant ultrafilters. In what is left 
we show that this is indeed possible. 
As a warm up exercise we begin describing how to turn a general $\kappa$-complete ultrafilter into a Galvin one by using a generalization of \emph{Mathias forcing}. The classical Mathias forcing dealing with subsets of $\omega$ appeared in \cite{MathiasHappy}, while the version that we will use here follows the template of \cite[Definition 3.1]{GartiShelah964}:
\begin{definition}[Generalized Mathias forcing]\label{PropertiesMathias}
\label{defgeneralmathias}
Let $\kappa$ be a regular cardinal 
and $\mathscr{U}$ a non-principal $\kappa$-complete filter over $\kappa$.

The forcing notion $\mathbb{M}_{\mathscr{U}}$ consists of pairs $(a,A)$ such that $a \in [\kappa]^{<\kappa}, A \in \mathscr{U}$ and $\sup(a)<\min(A)$.
For the order, one writes $(a_0, A_0) \leq (a_1, A_1)$ if and only if $a_0 \subseteq a_1, A_0 \supseteq A_1$ and $a_1\setminus a_0 \subseteq A_0$.
\end{definition}
The next is a brief account of the main properties of $\mathbb{M}_{\mathscr{U}}$:
\begin{proposition}[Properties of $\mathbb{M}_{\mathscr{U}}$]\label{PropertiesOfMathias}\hfill
\begin{enumerate}
	\item $\mathbb{M}_{\mathscr{U}}$ is $\kappa^+$-centered, provided $\kappa^{<\kappa}=\kappa$;
	\item $\mathbb{M}_{\mathscr{U}}$ is $\kappa$-directed-closed;
	\item $\mathbb{M}_{\mathscr{U}}$ is countably parallel closed;\footnote{I.e., every two decreasing sequences of conditions $\langle p_n\mid n<\omega\rangle$, $\langle q_n\mid n<\omega\rangle$ with $p_n\parallel q_n$ admit an upper bound. See \cite{AFramework}.}
	\item If $G\s \mathbb{M}_{\mathscr{U}}$ is a $V$-generic filter then $V[G]=V[a_G]$, where
	$$a_G:=\bigcup\{a\mid \exists A\in\mathscr{U}\,((a,A)\in G)\};$$
	\item The set $a_G$ \emph{diagonalizes} $\mathscr{U}$: i.e., $a_G\s^* A$  for every $A\in \mathscr{U}$. 
	\end{enumerate}
	In particular, if $\mathscr{U}$ is an ultrafilter then either $a_G\s^* A$ or 
$a_G\s^* \kappa\setminus A$ for all  $A\in \mathcal{P}(\kappa)^V$. \end{proposition}
The next proposition describes how to turn a $\kappa$-complete filter into one satisfying Galvin's property by means of $\mathbb{M}_\mathscr{U}$:
\begin{proposition}\label{MathiasandBase}
Let $\kappa$ be a regular  cardinal, $\mathscr{U}$ a $\kappa$-complete filter  over $\kappa$ and $G\s \mathbb{M}_\mathscr{U}$ a generic filter. Then the following hold in $V[G]$:
\begin{enumerate}
	\item $\mathscr{U}^*:=\{A\subseteq\kappa \mid a_G\subseteq^* A\}$ is a $\kappa$-complete filter;
	\item $\mathscr{U}\s \mathscr{U}^*$;
	\item $\mathrm{Gal}(\mathscr{U}^*,\kappa,\kappa^+)$ holds.
\end{enumerate}
\end{proposition}
\begin{proof}
 $\mathscr{U}^*$ is clearly a filter and $\mathscr{U}\s \mathscr{U}^*$ by virtue of Proposition~\ref{PropertiesMathias}(4). The argument for $\kappa$-completeness of $\mathscr{U}^*$ is essentially the same as the one for Clause~(3): Let $\l A_\alpha \mid \alpha <\kappa^+\r\s  \mathscr{U}^*$ and find $I\in [\kappa^+]^{\kappa^+}$ and $\alpha^*<\kappa$  such that $a_G\setminus \alpha^*\s A_\alpha$ for all $\alpha\in I$. Thus, $\bigcap_{\alpha\in I} A_\alpha\in\mathscr{U}^*$. 
\end{proof}
The above argument repeats the one from \cite[Proposition 4.5]{bgp}: the point is that $\mathbb{M}_\mathscr{U}$ creates a \emph{generating sequence} of length $1$. 
\begin{definition} 
	A family $\mathcal{A}=\langle x_\alpha\mid \alpha<\lambda\rangle\s \mathscr{U}$ is a \emph{generating sequence for $\mathscr{U}$} if for every $A\in \mathscr{U}$ there is $\alpha<\lambda$ such that $x_\alpha\s^* A.$ In addition, $\mathcal{A}$ is called \emph{strong generating} if it is $\subseteq^*$-decreasing.
\end{definition}

As demonstrated in \cite[\S4]{bgp} the analysis of (strong) generating sequences provides an effective way to produce certain Galvin-like configurations. The main obstacle, however, is to ensure that the departing $\kappa$-complete ultrafilter $\mathscr{U}$ extends to yet another $\kappa$-complete ultrafilter. This will be eventually addressed in Theorem~\ref{GitikShelahTheorem}.

\smallskip

Our next goal will be to iterate Mathias forcing over a given filter (and its extensions along the way) so that it will generate a $\kappa$-complete ultrafilter with a strong generating sequence of arbitrary length. 
This idea traces back to 
Kunen who employed it to separate the ultrafilter number $\mathfrak{u}$ from $2^{\omega}$ (see \cite[Ch. VII Question (A10)]{Kunen1980}). A similar argument, yet 
involving a more complex iteration, was considered in \cite{BrookeTaylor2017CardinalCA}. There 
the authors separate 
$\mathfrak{u}(\kappa)$ and $2^{\kappa}$ in a context where $\kappa$ is supercompact. 

 The \emph{na\"{i}ve} approach would consist of 
iterating Mathias forcing over and over with $\kappa$-complete filters. Unfortunately, this strategy is doomed to failure and so an 
additional structure on the forcing is required. Let us illustrate where the problem arises.  Suppose that $x_0$ is a Mathias set for a $\kappa$-complete filter  $\mathscr{U}$. Working in the generic extension $V[x_0]$ let $\mathscr{U}_0$ be a $\kappa$-complete filter extending $\{x_0\}\cup\mathscr{U}$ (e.g., 
by Proposition~\ref{MathiasandBase} we can take $\{X\in \mathcal{P}(\kappa)^{V[x_0]}\mid x_0\subseteq^*X\}$). Next, over $V[x_0]$, force a Mathias set  $x_1$ through $\mathscr{U}_0$ and  
working over the resulting extension $V[x_0,x_1]$ let $\mathscr{U}_1$ a 
$\kappa$-complete filter  extending $\{x_1\}\cup\mathscr{U}_0$. One can proceed in this fashion $\omega$-many times. Formally speaking, this is forced by the following 
full-support iteration $\l \mathbb{P}_n,\lusim{\mathbb{Q}}_n\mid n<\omega\r$: 
for each $n\geq 1$, $\lusim{\mathbb{Q}}_n$ is a $\mathbb{P}_n$-name for $\mathbb{M}_{\lusim{\mathscr{U}}_{n}}$ where $\lusim{\mathscr{U}}_n$ is a $\mathbb{P}_n$-name for a $\kappa$-complete ultrafilter  extending $\{\lusim{x}_{n}\}\cup\lusim{\mathscr{U}}_{n-1}$.

An essential obstacle arises at stage $\omega+1$. Here 
one needs to find a $\kappa$-complete filter which includes all the Mathias sets $\langle x_n\mid n<\omega\rangle$  constructed so far. However, notice that $\mathscr{W}:=\{X\subseteq\omega\mid \exists n<\omega\, (x_n\subseteq^* X)\}$ (i.e., the filter generated by the Mathias sets)  is not $\sigma$-complete: for if $\bigcap_{n<\omega}x_n\in\mathscr{W}$ then 
there would be some $n^*<\omega$ such that $x_{n^*}\subseteq^* \bigcap_{n<\omega} x_n$, 
 hence $x_{n^*}$ would be $\s^*$-included in $x_{n^*+1}$. This latter is certainly impossible in that $x_{n^*+1}$ is a Mathias set for a filter including $x_{n^*}$. 

\smallskip

For the moment, and as a warm up for Theorem~\ref{GitikShelahTheorem},  we show how to produce $\kappa$-complete filters with arbitrarily long strong generating sequences using  Mathias forcing.  

\begin{theorem}\label{ExtendFilters}
\label{propsupercom} Let $\mathscr{U}$ be a $\kappa$-complete filter over a Mahlo cardinal $\kappa$.

Then, for every $\lambda\in \mathrm{Ord}$ there is a $\kappa$-directed-closed and $\kappa^+$-cc poset $\mathbb{P}(\lambda)$ forcing 
that  $\mathscr{U}$ can be extended to a $\kappa$-complete filter $\mathscr{U}^*$ with a strong generating sequence $\langle x_\alpha\mid \alpha<\lambda\rangle$.
\end{theorem}
\begin{proof}
	Let $\mathscr{U}$ and $\lambda$ be as above.  Define a ${<}\kappa$-supported iteration $\mathbb{P}(\lambda)$, $\langle\mathbb{P}_\alpha,\lusim{\mathbb{Q}}_\beta\mid \beta<\alpha\leq\lambda\rangle$, as follows.
  Suppose that $\mathbb{P}_\alpha$ is defined for $\alpha<\lambda$. In $V^{\mathbb{P}_\alpha}$ we will define a filter $\mathscr{U}^*_\alpha$ over $\kappa$  which we will prove to be $\kappa-$complete. Bearing this in mind, we shall let $\lusim{\mathbb{Q}}_\alpha$ be a $\mathbb{P}_\alpha$-name for $\mathbb{M}_{\mathscr{U}^*_\alpha}$ and denote by $x_\alpha:=a_{G_{\mathbb{Q}_\alpha}}$ the  generic Mathias set added after forcing with $\lusim{\mathbb{Q}}_\alpha$. 
  
  For $\alpha=0$, we let $\mathscr{U}_0:=\mathscr{U}$. At successor $\alpha+1$, in $V^{\mathbb{P}_{\alpha+1}}$ we have $x_\alpha$, then we let $\lusim{\mathscr{U}}_{\alpha+1}$ be the $\mathbb{P}_{\alpha+1}-$name for the filter generated by $\lusim{x}_\alpha$.  By proposition \ref{MathiasandBase}, we have that $0_{\mathbb{P}_{\alpha+1}}\Vdash\lusim{\mathscr{U}}_\alpha\subseteq\lusim{\mathscr{U}}_{\alpha+1}$ and $\lusim{\mathscr{U}}_{\alpha+1}$ is $\kappa-$complete. As for the limit stages, let us split into cases
  \begin{claim}
  Suppose that $\alpha<\lambda$ is such that  $\cf(\alpha)\geq\kappa$. Then $\langle x_\beta\mid \beta<\alpha\rangle$  generates a $\kappa$-complete filter $\mathscr{U}^*_\alpha$ in $V^{\mathbb{P}_\alpha}$ that extends $\mathscr{U}^*_{\beta}$ for every $\beta<\alpha$.
  \end{claim}
  \begin{proof}
  Let $\{X_i\mid i<\mu<\kappa\}\subseteq \mathscr{U}^*_\alpha$. For each $i<\mu$ there is $\beta_i<\alpha$ such that $x_{\beta_i}\subseteq^* X_i$. Since the cofinality of $\alpha$ is at least $\kappa$,  $\sup_{i<\mu}\beta_i=:\beta^*<\alpha$. Hence, $x_{\beta^*}\subseteq^* x_{\beta_i}\subseteq^* X_i$ for every $i<\mu$. It follows that for some $\epsilon_i<\kappa$, $x_{\beta^*}\setminus\epsilon_i\subseteq X_i$. Take $\epsilon^*=\sup_{i<\mu}\epsilon_i<\kappa$, then $x_{\beta^*}\setminus\epsilon^*\subseteq\cap_{i<\mu}X_i$, by definition, $\cap_{i<\mu}X_i\in \mathscr{U}^*_\alpha$. Also, for every $\beta<\alpha$, $\mathscr{U}^*_{\beta}\subseteq\mathscr{U}^*_{\beta+1}$ and $\mathscr{U}^*_{\beta+1}$ is by definition the filter generated $x_\beta$ which is clearly a subset of $\mathscr{U}^*_\alpha$.
  \end{proof}
  \begin{claim}
  Suppose that $\alpha<\lambda$ is limit such that  $cf(\alpha)<\kappa$, and let $$\alpha=\kappa^{\delta_1}\gamma_1+...+\kappa^{\delta_n}\gamma_n$$ be the Cantor normal form of $\alpha$. Consider the following cofinal subset of $\alpha$: $I_\alpha:=\{\kappa^{\delta_1}\gamma_1+...+\kappa^{\delta_n}\gamma\mid \gamma<\gamma_m\}$. Then  $x:=\cap_{i\in I_\alpha}x_i\in V^{\mathbb{P}_\alpha}$ is unbounded in $\kappa$. In particular, $x\subseteq^* x_\beta$ for every $\beta<\alpha$, and the filter $\mathscr{U}^*_\alpha$ generated by $x$, is a $\kappa$-complete filter in $V^{\mathbb{P}_\alpha}$ which extends $\mathscr{U}^*_\beta$ for every $\beta<\alpha$.
  \end{claim}
  \begin{proof}
  Let $p\in\mathbb{P}_\alpha$ and $\alpha_0<\kappa$. We shall now proceed with a density argument to prove that there are $\alpha_0<\gamma^*<\kappa$  and $p\leq p_{fin}$ such that $p_{fin}\Vdash \gamma^*\in\cap_{i\in I_{\alpha}}\lusim{x}_i$. Construct two sequences $\l M_\rho\mid \rho<\kappa\r$ and $\l q_\rho\mid \rho<\kappa\r$ such that:
  \begin{enumerate}
      \item $M_\rho\prec H(\theta)$ for $\theta$ regular and sufficiently large.
      \item $M_{\rho}$ is increasing and continuous.
      \item $|M_\rho|:=\gamma_\rho<\kappa$. $\gamma_0$ is regular and also $\gamma_{\rho+1}$.
      \item If $\gamma_\rho$ is regular then $M_\rho^{<{\gamma_\rho}}\subseteq M_\rho$.
      \item $\alpha,p,\mathbb{P}_\alpha,\kappa\in M_\rho$.
      \item If $\beta\in M_{\rho}\cap \alpha+1$ has cofinality less than $\kappa$, then $cf(\beta)\cup I_{\beta}\subseteq M_{\rho+1}$.
  \end{enumerate}
  For the condition $q_\alpha$ we require that:
  \begin{enumerate}
      \item $q_\rho$ is increasing and continuous.\footnote{i.e. for limit $\rho$, we let $q_\rho(\gamma)=\l \cup_{\rho'<\rho}\lusim{a}^{q_{\rho'}}_\gamma,\cap_{\rho'<\rho}\lusim{A}^{q_{\rho'}}_\gamma\r$.}
      \item $q_\rho$ is $M_\rho$-generic for $\mathbb{P}_\alpha$, namely, for every dense open $D\subseteq\mathbb{P}_\alpha$, $D\in M_\rho$, $q_{\rho}\in D$. 
     \item $supp(q_{\rho})=M_\rho\cap\alpha$.
     \item $q_{\rho}\in M_{\rho+1}$
      
  \end{enumerate}
  For this construction we need that $\kappa$ is Mahlo. Such  condition exists by $\kappa-$closure of $\mathbb{P}_\alpha$ and by standard construction of an increasing sequence of conditions in $M$. Let $\mu^*<\kappa$ be regular such that $|M_{\mu^*}|=\mu^*$. Denote $M^*=M_{\mu^*}$ and $p^*=q_{\mu^*}=\sup_{i<\mu^*} q_{i}$.
  \begin{claim} For every $\beta\in Supp(p^*)$ the following hold:
 \begin{enumerate}
 \item  There is $\delta_{\beta}$ such that $p^*\restriction \beta\Vdash \lusim{a}^{p^*}_\beta\subseteq\delta_\beta$.
 \item If $\beta=\beta_0+1$ is successor, then there is $\epsilon_\beta$ such that $p^*\restriction \beta\Vdash \lusim{x}_{\beta_0}\setminus \epsilon_\beta\subseteq \lusim{A}^{p^*}_\beta$. 
 \item If $\beta\in Supp(p^*)$ is of cofinality less than $\kappa$ then there is $\epsilon_{\beta}$ such that $p^*\restriction \beta\Vdash \cap_{i\in I_\beta}\lusim{x}_i\setminus \epsilon_\beta\subseteq \lusim{A}^{p^*}_\beta$.
     \item If $\beta\in Supp(p^*)$ is of cofinality at least $\kappa$, then there $i_{\beta}<\beta$ 
     and $\epsilon_\beta$ such that $p^*\restriction \beta\Vdash \lusim{x}_{i_{\beta}}\setminus\epsilon_\beta\subseteq\lusim{A}^{p^*}_\beta$.
 \end{enumerate}
  \end{claim}
  \begin{proof}
  To see $(1)$, since $\beta\in Supp(p^*)=M^*\cap\alpha$, find $\xi_0<\mu^*$ such that $\beta\in M_{\xi_0}\cap\alpha$. In $M_{\xi_0+1}$ we can define the dense open set
  $$D_{\xi_0}:=\{q\in\mathbb{P}_{\alpha}\mid \exists \delta<\kappa. q\restriction \beta\Vdash \lusim{a}^{q_{\xi_0}}_\beta\subseteq\delta\}$$
  Since $q_{\xi_0+1}$ is $M_{\xi_0}-$generic, $q_{\xi_0+1}\in D_{\xi_0}$. For every $\xi_0\leq i<\mu^*$ we have that $q_i\in M_{i+1}$, hence we can define in $M_{i+1}$ the dense open set $$D_{i+1}:=\{q\in\mathbb{P}_{\alpha}\mid \exists \delta<\kappa. q\restriction \beta\Vdash \lusim{a}^{q_{i}}_\beta\subseteq\delta\}$$ By genericity, $q_{i+1}\in D_{i+1}$. For every such  $i$, pick $\delta^{(i)}$ witnessing $q_{i+1}\in D_{i+1}$. Let $\delta_\beta=\sup_{i<\mu^*}\delta^{(i)}<\kappa$, by continuity, $p^*\restriction \beta\Vdash \lusim{a}^{p^*}_\beta\subseteq \delta_\beta$.
  
 The proof of $(3),(4)$ is similar to $(1)$. Just note that by definition of $\mathscr{U}_\beta$, it is the filter generated by $\lusim{x}_{\beta_0}$ and replace the dense $D_{i+1}$ by
 $$E_{i+1}:=\{q\in\mathbb{P}_\alpha\mid \exists\epsilon. q\restriction \beta\Vdash \lusim{x}_{\beta_0}\setminus\epsilon\subseteq \lusim{A}^{q_i}_\beta\}$$
 Finally, to see $(4)$, follow a similar path by defining 
 $$F_{i+1}:=\{q\in\mathbb{P}_\alpha\mid \exists j_\beta<\beta.\exists\epsilon. q\restriction \beta\Vdash \lusim{x}_{j_\beta}\setminus\epsilon\subseteq \lusim{A}^{q_i}_\beta\}$$
 Pick for every $i<\mu^*$, $j_{\beta,i}\in M_{\rho^*}\cap\beta$ witnessing $q^*\in F_{i+1}$ and since the cofinality of $\beta$ is at least $\kappa$ we can take the sup to find a single $i_{\beta}$. Now choose the epsilons as before. Since $j_{\beta,i}\in M^*$ and unbounded in $i_{\beta}$,  there is $i_0<\mu^*$ such that for every $i_0\leq i<\mu^*$, the cantor normal form of  $j_{\beta,i}$ is a continuation of the on of $i_\beta$. Hence the $\delta_i'$s belong to $M^*$ 
  \end{proof}
  By $(1)$, for every $\beta\in Supp(p^*)$ we have $\delta_{\beta}<\kappa$, pick $\delta^*=\sup\delta_\beta$. By $(2),(3)$ and $(4)$ 
  we choose $\epsilon_\beta$ and set $\epsilon^*=\sup\epsilon_\beta<\kappa$. Pick $\gamma^*\in A^{p^*}_0$ above $\epsilon^*,\delta^*,\alpha_0$ and define $p_{fin}$. Define the support of $p_{fin}$ to be $Supp(p^*)\cup\{i_\beta\mid \beta\in Supp(p^*),cf(\beta)\geq\kappa\}$. Define
  $$p_{fin}(\gamma):=\begin{cases}\l\lusim{a}^{p^*}_\gamma\cup\{\gamma^*\},\lusim{A}^{p^*}_\gamma\setminus\gamma^*\r, & \gamma\in Supp(p^*),\\ \l \{\gamma^*\},\kappa\setminus\gamma^*+1\r, &otherwise.\end{cases}$$
  Clearly $p_{fin}\Vdash \gamma^*\in \cap_{i\in I_\alpha}\lusim{x}_i$. It remains to argue that $p_{fin}$ is an extension of $p^*$:
  Indeed $p_{fin}(0)\geq p^*_0$ since $\gamma^*\in A_0$ and $\gamma^*>\delta_0>\sup(a^{p^*}_0)$. Suppose that $\beta\in Supp(p^*)$ and that $p_{fin}\restriction \beta\geq p^*\restriction \beta$. If $\beta=\beta_0+1$ is successor then $p_{fin}\restriction\beta\Vdash \gamma^*\in \lusim{x}_{\beta_0}\setminus\epsilon^*\subseteq \lusim{A}^{p^*}_\beta$, hence $p_{fin}\restriction\beta\Vdash p_{fin}(\beta)\geq p^*(\beta)$. If $\beta$ is limit of cofinality less than $\kappa$, then since $\beta\in M^*$, we have that $I_{\beta}\subseteq M^*\cap\beta$, hence by induction
  $p_{fin}\restriction\beta\Vdash \gamma^*\in\cap_{i\in I_{\beta}}\lusim{x}_{\beta}\setminus\epsilon^*\setminus\lusim{A}^{p^*}_\beta$. Finally, if the cofinality of $\beta$ is at least $\kappa$, then there is $i_\beta<\kappa$ such that $p^*\restriction\beta\Vdash \gamma^*\in x_{i_\beta}^*\setminus\epsilon^*\subseteq \lusim{A}^{p^*}_\beta$.
  \end{proof}
This completes the proof of the theorem.  
\end{proof}
Let us briefly describe a natural, yet unfruitful,  strategy to make $\mathscr{U}^*$ become a $\kappa$-complete ultrafilter. 
Given a Laver-indestructible supercompact cardinal $\kappa$ and $\alpha<\lambda$ 
force over $V^{\mathbb{P}_\alpha}$ with $\mathbb{M}_{\mathscr{V}_\alpha}$, where $\mathscr{V}_\alpha$ is an extension of $\mathscr{U}^*_{\alpha}$ (in $V^{\mathbb{P}_\alpha}$) to a $\kappa$-complete ultrafilter. Note that since $\mathbb{P}_\alpha$ is $\kappa$-directed-closed, $\kappa$ is still supercompact in the corresponding extension and thus the choice of $\mathscr{V}_\alpha$ is available. 
In addition, if the length of the iteration is $\lambda=\kappa^+$ then the union of the  (tower of) $\kappa$-complete ultrafilters generated along the way will be also an ultrafilter, $\mathscr{U}_\infty$. The problem with this approach is that we lose control upon $\kappa$-completeness of $\mathscr{U}_\infty$. Indeed, even the union of the first $\omega$-many ultrafilters generated 
 might not be $\kappa$-complete, as we argued in the discussion preceding Theorem~\ref{ExtendFilters}.
 
 \smallskip

 The approach of Theorem~\ref{GitikShelahTheorem} is to iterate $\mathbb{M}_\mathscr{U}$ more carefully so that we have complete control on the completeness of the  ultrafilter $\mathscr{U}_\infty$.

 \begin{definition}
 	For  $f\colon \kappa\rightarrow\kappa$ and an ultrafilter $\mathscr{U}$ over $\kappa$ we say that  \emph{$f$ is constant  $\mathrm{mod}(\mathscr{U})$} if there is $\gamma<\kappa$ such that $f^{-1}\{\gamma\}\in\mathscr{U}$. Similarly, 
 \emph{$f$ is $1$-$1\mathrm{mod}(\mathscr{U})$} if there is $X\in \mathscr{U}$ such that $|f^{-1}\{\gamma\}\cap X|<\kappa$ for every $\gamma<\kappa$. 
 \end{definition}

\begin{definition}
A $\kappa$-complete ultrafilter $\mathscr{U}$ is called a \emph{$P$-point} if every function $f:\kappa\rightarrow\kappa$ that is not constant $\mathrm{mod}(\mathscr{U})$  is 
$1$-$1\mathrm{mod}(\mathscr{U})$.
\end{definition}
\begin{remark}
	$\mathscr{U}$ is a $P$-point if and only if every sequence  $\langle X_\alpha\mid \alpha<\kappa\rangle$ of elements in $\mathscr{U}$ has a pseudo intersection in $\mathscr{U}$; namely, there is $X\in \mathscr{U}$ such that $X\subseteq^* X_\alpha$ for every $\alpha<\kappa$.
\end{remark}

    
\begin{lemma}\label{P point is Galvin}
Every $\kappa$-complete ultrafilter $\mathscr{U}$ with a  generating sequence of size $\kappa^+$ is a $P$-point. In particular,  $\mathrm{Gal}(\mathscr{U},\kappa,\kappa^+)$.
\end{lemma}
\begin{proof}
 Let $\mathcal{A}=\l A_\alpha\mid\alpha<\kappa^+\r$ be a generating sequence of $\mathscr{U}$. To see it is $P$-point suppose that $\l X_\alpha\mid \alpha<\kappa\r$ is any $\kappa-$sequence of members of $\mathscr{U}$. By definition of generating sequence, for each $\alpha<\kappa$ there is $\beta_\alpha<\kappa^+$ such that $A_{\beta_\alpha}\subseteq^* X_\alpha$. Consider $\beta^*=\sup_{\alpha<\kappa}\beta_\alpha<\kappa^+$. Then $A_{\beta^*}$ is a pseudo intersection of the sequence $\l X_\alpha\mid \alpha<\kappa\r$. For the in particular claim use the fact that every $P$-point ultrafilter $\mathscr{U}$ satisfies $\mathrm{Gal}(\mathscr{U},\kappa,\kappa^+)$ (see \cite[Proposition~5.13]{MR4393795}).
 \end{proof}
In analogy to Theorem \ref{MakingEverythingNormal}, next we show that every ground model $\kappa$-complete ultrafilter can be extended to a \emph{well-behaved} one: to wit, to a $P$-point. By virtue of the above lemma, this gives an alternative (and less severe way) to trasnmute an arbitrary $\kappa$-complete ultrafilter into a Galvin one. 
Unlike Theorem~\ref{MakingEverythingNormal},  the models we produce this time have the extra feature that $2^\kappa$ can be made arbitrarily large. Our construction owes much to previous work of Gitik and Shelah \cite{MR1632081}. 
Recall that $\kappa$ is \emph{almost huge with target $\lambda$} if there is $j\colon V\rightarrow M$ such that $\crit(j)=\kappa$, $j(\kappa)=\lambda$ and $M^{<\lambda}\s M$. 
\begin{theorem}\label{GitikShelahTheorem} 
Assume that the $\mathsf{GCH}$ holds and  suppose that $\kappa$ is an \linebreak almost huge cardinal with measurable target $\lambda$. Then for every $\delta<\lambda$ there is a forcing extension $V[G_\delta]$ where 
$2^\kappa=\delta$ and every ground model $\kappa$-complete ultrafilter extends to a $P$-point ultrafilter in $V[G_\delta]$. In addition, $V[G_\delta]_\lambda$ models the same configuration and $\kappa$ is supercompact there.
\end{theorem}
\begin{proof}
 Let $j\colon V\rightarrow M$ be such that $\crit(j)=\kappa$, $M^{<\lambda}\s M$ and $j(\kappa)=\lambda$. Recall that $\lambda$ is  assumed to be measurable in the ground model, hence we can let $\mathcal{U}$ be a measure on $\lambda$.  Let us define an Easton support iteration $\l \mathbb{P}_\alpha, \lusim{\mathbb{Q}}_\beta\mid \beta\leq \kappa,\, \alpha\leq \kappa+1\r$, where  $\lusim{\mathbb{Q}}_\alpha$ is trivial unless $\alpha$ is measurable in $V^{\mathbb{P}_\alpha}$, in which case $\lusim{\mathbb{Q}}_\alpha$ is a $\mathbb{P}_\alpha$-name for the two-step  iteration $\lusim{\mathbb{Q}}_{\alpha,0}\ast\lusim{\mathbb{Q}}_{\alpha,1}$ defined as follows: $\lusim{\mathbb{Q}}_{\alpha,0}$ is the atomic forcing choosing some ordinal $F(\alpha)<\kappa$ followed by $\lusim{\mathbb{Q}}_{\alpha,0}$, an ${<}\alpha$-supported iteration $\langle \mathbb{R}^\alpha_{\beta}, \mathbb{S}^\alpha_\gamma\mid \beta\leq F(\alpha),\, \gamma<F(\alpha)\rangle$ defined as follows:  At  each step $\beta\leq F(\alpha)$, $\mathbb{S}^\alpha_\beta$ is trivial unless $\Vdash_{\mathbb{P}_\alpha\ast\mathbb{R}^\alpha_{\beta}}``\alpha$ is measurable'', 
 in which case $\mathbb{S}^\alpha_\beta$ is forced to be the ${<}\alpha$-supported product $\prod_{\lusim{\mathscr{V}}}\mathbb{M}_{\lusim{\mathscr{V}}}$, where $\lusim{\mathscr{V}}$ ranges over all $\mathbb{P}_\alpha\ast\mathbb{R}^\alpha_{\beta}$-names for an $\alpha$-complete ultrafilter over $\alpha.$ 
 
 Since $\mathbb{Q}_{\alpha,1}$ is an ${<}\alpha$-supported iteration of $\alpha^+$-stationary c.c., ${<}\alpha$-closed and {countably-parallel closed} 
 forcing (cf.~Proposition~\ref{PropertiesOfMathias}) it follows from \cite[Theorem~1.2]{AFramework} that  $\mathbb{Q}_{\alpha,1}$ is $\alpha^+$-cc. 
 Also, $\mathbb{P}_\kappa$ is $\kappa$-cc because $\kappa$ is Mahlo and $\mathbb{P}_\kappa$ is Easton-supported. Let $G_\kappa\subseteq \mathbb{P}_\kappa$ be a $V$-generic filter. Then, due to closure of $M$ under ${<}\lambda$-sequences and $\kappa$-ccness of $\mathbb{P}_\kappa$, $V[G_\kappa]$ and $M[G_\kappa]$ agree up to $\lambda$ (see, e.g., \cite[Proposition 8.4]{MR2768691}). Consider $j(\mathbb{P}_\kappa):=\mathbb{P}'_{j(\kappa)}$.
 \begin{claim}
 For every $\rho<\lambda$, $M[G_\kappa\ast \{\rho\}]^{\mathbb{Q}_{1,\kappa}}\models ``\kappa$ is measurable''. 
 
 In fact, $\kappa$ is ${<}\lambda$-supercompact in $M[G_\kappa\ast\{\rho\}]^{\mathbb{Q}_{1,\kappa}}$, hence also in \linebreak $V[G_\kappa\ast\{\rho\}]^{\mathbb{Q}_{1,\kappa}}$, and thus $\kappa$ is fully supercompact in  $(V[G_\kappa\ast\{\rho\}]^{\mathbb{Q}_{1,\kappa}})_\lambda.$
 \end{claim}
 \begin{proof}
  In the ground model $V$, let $\max(\rho,2^{2^{\kappa}})\leq \theta<\lambda$. Let $U$ be a fine normal measure over $P_\kappa(\theta)$ and let $j_U\colon V\rightarrow M$ be the corresponding elementary embedding. Then ${}^{\theta}M\subseteq M$. Put $\mathbb{P}'_{j_U(\kappa)}:=j_U(\mathbb{P}_\kappa)$ and $\mathbb{Q}'_{j_U(\kappa)}:=j_U(\mathbb{Q}_\kappa)$. Let $G(\mathbb{Q}_{1,\kappa})$ a $V[G_\kappa][\{\rho\}]$-generic, and let us lift $j_U$ to the model $V[G_{\kappa}\ast \{\rho\}\ast G(\mathbb{Q}_{1,\kappa})]$. Note that in $M[G_{\kappa}\ast \{\rho\}\ast G(\mathbb{Q}_{1,\kappa})]$ we have $2^\kappa\geq \rho$, since at each step of the iteration $\mathbb{Q}_{1,\kappa}$ we add a new subset to $\kappa$. Hence, by closure under $\rho$-sequences, the degree of closure of  the forcing $$\mathbb{P}'_{j_{U}(\kappa)}/[G_{\kappa}\ast\{\rho\}*G(\mathbb{Q}_{1,\kappa})]$$ is at least $\rho^{+}$, even in $V[G_{\kappa}\ast\{\rho\}\ast G(\mathbb{Q}_{1,,\kappa})]$. Also note that every dense set in $\mathbb{P'}_{j(\kappa)}/G_{\kappa}\ast \{\rho\}\ast G(\mathbb{Q}_{1,\kappa})$ is represented by a function $f:\mathcal{P}_\kappa(\rho)\rightarrow \mathcal{P}(\mathbb{P}_\kappa)$. Since there are in total $(2^{\kappa})^{\rho}=\rho^+$ of such functions we can construct a master sequence which induces an $M[G_\kappa\ast\{\rho\}\ast G(\mathbb{Q}_{1,\kappa})]$-generic filter $G_{<\lambda}$. Next, the top-most forcing $\mathbb{Q}_{0,\kappa}\ast \mathbb{Q}_{1,\kappa}$ is $\kappa^+$-cc and when its length is restricted to some $\rho$ its {size becomes  $\rho\cdot 2^{2^\kappa}$}. Hence every maximal antichain in $M[G_{<\lambda}\ast\{j(\rho)\}]$ for $j(\mathbb{Q}_{1,\kappa})$ is represented by a function $F\colon \mathcal{P}_{\kappa}(\rho)\rightarrow \mathcal{P}_{\kappa}((\mathbb{Q}_{1,\kappa})_{G_{\kappa}\ast\{\rho\}})\in V[G_{\kappa}\ast\{\rho\}]$ and since $(\mathbb{Q}_{1,\kappa})_{G_{\kappa}\ast\{\rho\}}$ is of size $\rho\cdot 2^{2^\kappa}=\theta$ there are $\theta^{\theta}=\theta^+$ many such functions. The remaining argument is standard.
 \end{proof}
 For each $\rho<\lambda$,  $(\mathbb{Q}_{1,\kappa})_{G_\kappa\ast \{\rho\}}$ is $\kappa^{+}$-cc and of size less than $\lambda$. Thus, there are less than $\lambda$-many nice names for subsets of $\kappa$. Denote these names by  $\l \lusim{A}^{\rho}_\tau\mid \tau<\theta_\rho\r$. We can find in $V[G_{\kappa}]$ an enumeration $\l \lusim{A}_{\tau}\mid \tau<\lambda\r$ of subsets of $\kappa$ such that for every $\tau_1\leq \tau_2$, there are $\delta_1,\delta_2$ such that $\lusim{A}_{\tau_1}=\lusim{A}^{\rho(\tau_1)}_{\delta_1}$ and $\lusim{A}_{\tau_2}=\lusim{A}^{\rho(\tau_2)}_{\delta_2}$, and $\rho(\tau_1)\leq\rho(\tau_2)$. For each $\tau<\lambda$ for which $\rho(\tau)<\lambda$ has been defined, let $C'$ be the club of closure points of $\rho(\tau)$. Since $\lambda$ is measurable there is $S\in \mathcal{U}$ concentrating on inaccessibles. Hence, there for $\mathcal{U}$-many $\delta\in S$ such that $\rho``\delta\s \delta$. In particular, the sequence $\l \lusim{A}_\tau\mid\tau<\delta\r$ 
 codes all the $(\mathbb{Q}_{1,\kappa})_{G_\kappa*\{\delta\}}$-names for subsets of $\kappa$.\footnote{Specifically, if $\sigma$ is a $(\mathbb{Q}_{1,\kappa})_{G_\kappa*\{\delta\}}$-name for a subset of $\kappa$ then there is $\tau<\delta$ such that $0\Vdash_{(\mathbb{Q}_{1,\kappa})_{G_\kappa*\{\delta\}}}\sigma=\lusim{A}_\tau.$} 
 
 \smallskip
 
Let $\epsilon<\lambda$ be an ordinal above the generators of all $\kappa$-complete ultrafilters. More precisely, for each $\kappa$-complete ultrafilter $U\in V$ over $\kappa$ let $\varepsilon_U\in \lambda\cap \bigcap j``U$ and define $\epsilon:=\sup_{U}\varepsilon_U$. Note that $\epsilon<\lambda$ as $\lambda$ is inaccessible in $V$.
 
 For every  $\delta'<\delta$ and every $(\mathbb{Q}_{1,\kappa})_{G_\kappa\ast\{\delta' \}}$-name (i.e., a $(\mathbb{Q}_{1,\kappa})_{G_\kappa\ast\{\delta\}}\restriction\delta'$-name)  $\lusim{U}$ for a $\kappa$-complete ultrafilter over $\kappa$, 
 let us define  $r_{\lusim{U}}\in \mathbb{M}_{j(\lusim{U})}$ as: 
 $$r_{\lusim{U}}=\l a_{\lusim{U}}\cup (A_{\lusim{U}}\cap \epsilon),{A}_{\lusim{U}}\setminus (\epsilon +1)\r,$$ where
 \begin{itemize}
     \item $a_{\lusim{U}}$ is the standard $(\mathbb{Q}_{1,\kappa})_{G_\kappa\ast\{\delta\}}$-name for the Mathias set for $\mathbb{M}_{\lusim{U}}$;
     \item $A_{\lusim{U}}$ is a name for the set $\bigcap j``\lusim{U}$.
 \end{itemize}
   Note that $r_{\lusim{U}}$ is a condition in $\mathbb{M}_{j(\lusim{U})}$: First, the trivial condition of $\mathbb{M}_{j(\lusim{U})}$ forces $``\min A_{\lusim{U}}=\kappa$'', hence $a_{\lusim{U}}\cup (A_{\lusim{U}}\cap \epsilon)$ is a legitimate stem. Second, this condition also forces $j(\lusim{U})$ to be $j(\kappa)$-complete, hence $\bigcap j``\lusim{U}\in j(\lusim{U})$.

  The key feture of $r_{\lusim{U}}$ is that $0_{j(\mathbb{M}_{\lusim{U}})}$ forces $a_{j(\lusim{U})}$ to contain all generators $\varepsilon_W$, $W\in V$, that are forced to belong to $A_{\lusim{U}}$. 
  
   \smallskip
Let $q_{\delta}$ be a condition in $\mathbb{Q}'_{j_U(\kappa)}$ with support $j``\delta$ (hence of cardinality ${<}j(\kappa)$) be the following condition: For every $\delta'<\delta$ and a $(\mathbb{Q}_{1,\kappa})_{G_\kappa\ast\{\delta' \}}$-name for a $\kappa$-complete ultrafilter $\lusim{U}$ as above, 
 $$q_{\delta}\restriction j(\delta')\Vdash q_{\delta}(j(\delta'))(j(\lusim{U}))=r_{\lusim{U}};$$
 for other coordinates (i.e., names for \emph{ghost ultrafilters}) $\lusim{W}$ we require that  $$q_{\delta}\restriction j(\delta')\Vdash \text{$``q_{\delta}(j(\delta'))(\lusim{W})$ is the trivial codition in $\mathbb{M}_{\lusim{W}}$''}.$$
 
 \smallskip
 
 A moment's reflection makes clear that $q_\delta$ is a \emph{master condition} for $G(\mathbb{Q}_{1,\kappa})$: namely,  $j(p)\leq q_\delta$ for every $p\in G({\mathbb{Q}_{1,\kappa}}).$ 
 In addition, the conditions $q_\delta$ are defined in a coherent way: namely, if $\rho\leq \delta$ are both in $C$ then $q_\delta\restriction\rho=q_\rho$.


\smallskip

Fix $\rho \in C\setminus(\epsilon+1)$. For every $\tau,\zeta<\rho$ 
let $\lusim{D}_{\tau,\zeta}\in M[G_{\kappa}\ast\{\rho\}]$ be a  $(\mathbb{Q}_{1,\kappa})_{G_\kappa*\{\rho\}}$-name  for the following dense open set  $$\{p\in\mathbb{P}'_{(\kappa,j_U(\kappa)]} \mid \exists s_\kappa\in\mathbb{Q}_{1,\kappa}\,\exists i\in 2\, (\l s_\kappa,p\r \Vdash^i_{\mathbb{P}'_{[\kappa,j_U(\kappa)]}}  \zeta\in j(\lusim{A}_\tau))\}.$$
    Since the trivial condition of $\mathbb{Q}_{1,\kappa}$ forces $\mathbb{P}'_{(\kappa,j_U(\kappa)]}$ to be $\rho^+$-closed it also forces that $\bigcap_{\tau,\zeta}\lusim{D}_{\tau,\zeta}$ is a name for a dense open set. In particular, there is some $(\mathbb{Q}_{1,\kappa})_{G_{\kappa}\ast\{\rho\}}$-name $p_{\rho}$ for a condition in $\bigcap_{\tau,\zeta}\lusim{D}_{\tau,\zeta}$ such that {$p_{\rho}\restriction j(\kappa)\Vdash p_{\rho}(j(\kappa))\geq q_\rho$.} Notice that
    $p_\rho$ has the property that for every $\tau,\zeta<\rho$ there is $s\in G_{\kappa}*\{\rho\}$ and $s_\kappa\in(\mathbb{Q}_{1,\kappa})_{G_\kappa\ast\{\rho\}}$ such that  $\langle s,s_\kappa,p_\rho\rangle||_{\mathbb{P'}_{[\kappa,j_U(\kappa))}}\zeta\in j(\lusim{A}_\tau)$.
    

  For each $(s,s_\kappa)\in G_\kappa \ast \lusim{\mathbb{Q}}_{\kappa}$, ordinals $\zeta<\epsilon$ and $\tau<\lambda$, and $i\in 2$ define $$A^{i}_{(s,{s}_\kappa,\zeta,\tau)}:=\{\rho<\lambda\mid \langle s, s_\kappa, p_\rho\rangle\Vdash^i\zeta\in j(\lusim{A}_\tau)\}.$$  Denote by $A^{2}_{(s, {s}_\kappa,\zeta,\tau)}$ the complement of the union of the above two sets. For each such quadruple $(s,s_\kappa,\zeta, \tau)$, let $i_{(s,s_\kappa,\zeta, \tau)}\in 3$ be the unique index $i$ for which $A^{i}_{(s,{s}_\kappa,\zeta,\tau)}\in \mathcal{U}$, the $\lambda$-complete measure on $\lambda$. Now let $$A:=\{\rho<\lambda\mid (\langle s,s_\kappa\rangle \in G\ast \lusim{\mathbb{Q}}_\kappa\restriction\rho\,\wedge\, \max(\zeta,\tau)<\rho)\,\Rightarrow\, \rho\in A^{i_{(s,s_\kappa,\zeta,\tau)}}_{(s,s_\kappa,\zeta,\tau)}\}.$$
   We claim that $A\in\mathcal{U}$: In effect, for every $\langle s,s_\kappa\rangle \in G\ast j(\mathbb{Q}_\kappa)\restriction\lambda=G\ast \mathbb{Q}_\kappa$ and $\zeta,\tau<\lambda$,  $A^{i_{(s,s_\kappa,\zeta,\tau)}}_{(s,s_\kappa,\zeta,\tau)}\in\mathcal{U}$. Hence, $\lambda\in j_\mathcal{U}(A^i_{(s,s_\kappa,\zeta,\tau)})$, and thus $\lambda\in j_\mathcal{U}(A)$. 
   
   Put $C^*:=A\cap C$. Let $\rho,\rho'\in C^*$ with $\rho<\rho'$, $\langle s,s_\kappa\rangle \in G\ast (\lusim{\mathbb{Q}_{\kappa}}\restriction \rho)$ and $\zeta, \tau<\rho$. By definition of $A$,  $\rho,\rho'\in A^{i_{(s,s_\kappa,\zeta,\tau)}}_{(s,s_\kappa,\zeta,\tau)}$, hence  
    $$
    \text{$\langle s, s_\kappa, p_\rho\rangle \Vdash^i \zeta\in j(\lusim{A}_\tau)$ iff $\langle s, s_\kappa, p_{\rho'}\rangle \Vdash^i \zeta\in j(\lusim{A}_\tau)$,}
    $$
    and also
    $$
    \text{$\langle s, s_\kappa, p_\rho\rangle \Vdash \zeta\in j(\lusim{A}_\tau)$ iff $\langle s, s_\kappa, p_{\rho'}\rangle \Vdash \zeta\in j(\lusim{A}_\tau).$}
   $$
    
    Next, for all $\rho\in C^*$ and $\zeta<\epsilon$ consider
   $$\mathscr{U}_{\rho,\zeta}:=\{(\lusim{A}_\tau)_{G_{\kappa}\ast\{\rho\}\ast G(\mathbb{Q}_{\kappa,1})}\s \kappa\mid  \exists\langle s,s_\kappa\rangle \in G\ast \{\rho\}\ast G(\mathbb{Q}_{\kappa,1})\, \langle s, s_\kappa, 
    p_\rho \rangle \Vdash\zeta\in j(\lusim{A}_\tau)\}.$$
    Since $\langle \lusim{A}_\tau\mid \tau<\rho\rangle$ is an enumeration of the $(\mathbb{Q}_{\kappa,1})_{G_\kappa\ast \{\rho\}}$-names, it is not hard to show that $\mathscr{U}_{\rho,\zeta}$ is a $\kappa$-complete ultrafilter in  $V[G_\kappa\ast\{\rho\}\ast G(\mathbb{Q}_{\kappa,1})]$. 
    
    Also, for each $\zeta<\epsilon$, $\langle \mathscr{U}_{\rho,\zeta}\mid \rho\in C^*\rangle$ defines a tower of ultrafilters: Suppose $\rho<\rho'\in C^*$ and let $A\in \mathscr{U}_{\rho,\zeta}$. Then, there is a pair $\langle s,s_\kappa\rangle$ such that $\langle s,s_\kappa, p_\rho\rangle \Vdash\zeta\in j(\lusim{A}_\tau)$. By our definition of $C^*$ this is also true when replacing $p_\rho$ by $p_{\rho'}$. Thus, $(\lusim{A}_{\tau})_{G_{\kappa}\ast\{\rho\}\ast G(\mathbb{Q}_{\kappa,1})}\in\mathscr{U}_{\rho',\zeta}.$ 
    
    Let $\delta$ be the limit of some sequence $\langle \rho_\alpha\mid\alpha <\kappa^+\rangle\s C^*$. From our previous comments, $\langle \mathscr{U}_{\rho_\alpha,\zeta}\mid\alpha<\kappa^+\rangle$ defines a tower of measures. Now, define $\mathscr{V}_{\delta,\zeta}:=\bigcup_{\alpha<\kappa^+} \mathscr{U}_{\rho_\alpha,\zeta}$. Since $V[G_\kappa\ast \{\rho_\alpha\}\ast G(\mathbb{Q}_{\kappa,1})\restriction\rho_\alpha]$ is a submodel of $V[G_\kappa\ast \{\delta\}\ast G(\mathbb{Q}_{\kappa,1})]$ and  $(\mathbb{Q}_{1,\kappa})_{G_\kappa\ast \{\delta\}}$ is $\kappa^+$-cc it is immediate that $\mathscr{V}_{\delta,\zeta}$ is a $\kappa$-complete ultrafilter.
    
    \begin{claim} 
    	 $\mathscr{V}_{\delta,\zeta}$ admits a strong generating sequence of size  $\kappa^+$.
    \end{claim} 
    \begin{proof}
    For each $\alpha<\kappa^+$,  $\mathscr{U}_{\rho_\alpha,\zeta}\in V[G_\kappa\ast \{\delta\}\ast G(\mathbb{Q}_{\kappa,1})\restriction \rho_\alpha]$ hence the iteration $\mathbb{Q}_{\kappa,1}$ at stage $\rho_{\alpha}+1$ shoots a Mathias set $x_\alpha\s \kappa$ (over the model $V[G_\kappa\ast \{\delta\}\ast G(\mathbb{Q}_{\kappa,1})\restriction \rho_\alpha]$) for the measure $\mathscr{U}_{\rho_\alpha,\zeta}$. Namely, $x_\alpha$ is almost included in every $A\in \mathscr{U}_{\rho_\alpha,\zeta}$.

    	We claim that $\langle x_\alpha\mid \alpha<\kappa^+\rangle$ is the sought strong generating sequence. First, for each $A\in \mathscr{V}_{\delta,\zeta}$ there is $\alpha<\kappa^+$ such that $A\in\mathscr{U}_{\rho_\alpha,\zeta}$ and so $x_\alpha\s^* A$. Second, $\langle x_\alpha\mid \alpha<\kappa^+\rangle$ is $\s^*$-decreasing: Fix $\alpha<\beta<\kappa^+$. 
    	We would like to show that $x_\alpha\in \mathscr{U}_{\rho_\beta,\zeta}$. Recall that $x_{\alpha}=(a_{\lusim{\mathscr{U}}_{\rho_\alpha,\zeta}})_{G_{\kappa}\ast \{\rho_{\beta}\}\ast G(\mathbb{Q}_{\kappa,1})\restriction\rho_{\beta}}$, hence we should check that for some $\l s,s_\kappa\r\in G_{\kappa}\ast\{\rho_{\beta}\}\ast G(\mathbb{Q}_{\kappa,1})\restriction\rho_{\beta}$ we have that:
    	$$\l s,s_{\kappa},p_{\rho_{\beta}}\r\Vdash \zeta\in j(a_{\lusim{\mathscr{U}}_{\rho_\alpha,\zeta}})$$
    	By elementarity of $j$, it follows that $j(a_{\lusim{\mathscr{U}}_{\rho_\alpha,\zeta}})=a_{j(\lusim{\mathscr{U}}_{\rho_\alpha,\zeta})}$ is the canonical $\mathbb{P}'_{j(\kappa)}\ast\{j(\rho_{\beta})\}\ast \mathbb{Q}'_{j(\kappa)}$-name for the Mathias generic of $\mathbb{M}_{j(\mathscr{U}_{\rho_{\alpha},\zeta})}$. By the definition of $p_{\rho_{\beta}}$ (which extends $q_{\rho_{\beta}}$, and the definition of $r_{\lusim{\mathscr{U}}_{\rho_\alpha,\zeta}}$), $$\l 0,p_{\rho_{\beta}}\r\Vdash a_{j(\lusim{\mathscr{U}}_{\rho_\alpha,\zeta})}\cap(\epsilon+1)\supseteq \bigcap j``\lusim{\mathscr{U}}_{\rho_{\alpha},\zeta}\cap (\epsilon+1).$$  
    	Working in $M[G_{\kappa}\ast\{\rho_{\alpha}\}\ast G(\mathbb{Q}_{\kappa,1})]$, we have that for every $A\in \mathscr{U}_{\rho_\alpha,\zeta}$, there is a name $A_\tau$ for $A$ such that $p_{\rho_\alpha}\Vdash \zeta\in j(\lusim{A}_\tau)$. Hence $p_{\rho_\alpha}\Vdash \zeta\in \bigcap j``\lusim{\mathscr{U}}_{\rho_\alpha,\zeta}$. Hence there is $\l s,s_\kappa\r\in G_\kappa\ast\{\rho_\alpha\}\ast G(\mathbb{Q}_{\kappa,1})\restriction \rho_\alpha$ such that
    	$$\l s,s_\kappa,p_{\rho_\alpha}\r \Vdash \zeta\in j(a_{\lusim{\mathscr{U}}_{\rho_{\alpha},\zeta}})$$
    	Since $\rho_\alpha,\rho_\beta\in C^*$, this means that 
    	$\l s,s_\kappa,p_{\rho_{\beta}}\r\Vdash \zeta\in j(a_{\lusim{\mathscr{U}}_{\rho_{\alpha},\zeta}})$ which by definition implies that $x_\alpha\in\mathscr{U}_{\rho_{\beta},\zeta}$.
 \end{proof}

   \begin{claim}
   Every $\kappa$-complete ultrafilter $\mathscr{U}$ from the ground model is extended by $\mathscr{V}_{\delta,\zeta}$, for some $\zeta<\epsilon.$
   \end{claim}
   \begin{proof}[Proof of claim]
   Let $\zeta<\epsilon$ be such that $\mathscr{U}=\{X\s \kappa\mid \zeta\in j(X)\}$ . Clearly, $\mathscr{U}\s \mathscr{U}_{\rho_\alpha,\zeta}$ for all $\alpha<\kappa^+$, hence $\mathscr{U}\s \mathscr{V}_{\delta,\zeta}$.
   \end{proof}
    This completes the proof of the theorem.
\end{proof}

\begin{remark}
     Note that every measurable cardinal $\kappa$ always carries a $\kappa$-complete   ultrafilter which is not a $P$-point. To see this, take any $\kappa$-complete ultrafilter $\mathscr{U}$ over $\kappa$, and a bijection $\phi:[\kappa]^2\rightarrow \kappa$ and define  $\mathscr{W}:=\phi_*(\mathscr{U}\times\mathscr{U})$. One can check that $\mathscr{W}$ is a $\kappa$-complete non-$P$-point ultrafilter. As witnessed by $L[\mathscr{U}]$, it is consistent that every ultrafilter is a finite power of a normal one (hence of a $P$-point), and such ultrafilters are always Galvin (see \cite[Corollary 5.29]{MR4393795}). If $\kappa$ is $\kappa$-compact\footnote{A cardinal $\kappa\geq \aleph_1$ is called \emph{$\kappa$-compact} if every $\kappa$-complete filter over $\kappa$ extends to a $\kappa$-complete ultrafilter.} then there is  a $\kappa$-complete ultrafilter which is not a finite power of a $P$-point (see, e.g., \cite[\S3.9]{Kan1}). 
\end{remark}

\begin{question}
Is it consistent to have a measurable cardinal carrying a $\kappa$-complete ultrafilter $\mathscr{U}$ such that $\mathrm{Gal}(\mathscr{U},\kappa,\kappa^+)$ but it is not Rudin-Keisler equivalent to a finite power of $P$-points? 
\end{question}
\begin{question}
Is it $\mathsf{ZFC}$-provable that a supercompact cardinal always admits a $\kappa$-complete non Galvin ultrafilter?  
\end{question}

\subsection{Galvin's property in the choiceless context}

Another way to exa\-mine Galvin's property at very large cardinals is to consider relatively small cardinals in \textsf{ZF}.
A typical example is $\aleph_1$ under \textsf{AD}.
Indeed, Solovay proved that both $\aleph_1$ and $\aleph_2$ are measurable under \textsf{AD} and that $\mathscr{D}_{\aleph_1}$ is a normal ultrafilter (see \cite[Theorem~33.12]{Jech2003}).
In fact, $\aleph_1$ is $\aleph_2$-supercompact under \textsf{AD}, by a result of Martin (see \cite[p. 401]{Kan}). 
Moreover, under $\mathrm{AD}_\mathbb{R}$, $\aleph_1$ is $\gamma$-supercompact for all $\gamma<\Theta$ \cite[Theorem 28.22]{Kan}. 
However, the classical proof of Galvin's theorem employs the Axiom of Choice in a crucial way. 

In private communication \cite{Chan}, W. Chan informed us that  it is possible to get Galvin's property at many cardinals even in the absence of \textsf{AC}. 
The argument involves Steel's notion of \emph{boldface $\mathsf{GCH}$ cardinal} \cite[p. 1678]{MR2768698}:
\begin{definition}
A cardinal $\kappa$ is called \emph{boldface $\mathsf{GCH}$} if there is no injection $F\colon \kappa^+\rightarrow \mathcal{P}(\kappa)$.
\end{definition}
 In \cite[Theorem~8.26]{MR2768698} Steel shows that under $``V=L(\mathbb{R})+\mathsf{AD}$'' every cardinal $\kappa<\Theta$ is boldface \textsf{GCH.} Moreover, Chan indicated that using \emph{Moschovakis coding lemma} 
 one can show that \textsf{AD} alone implies that eve\-ry $\kappa<\Theta^{L(\mathbb{R})}$ is boldface \textsf{GCH}. 
\begin{proposition}\label{thmad}
If $\kappa$ is boldface $\mathsf{GCH}$ and $\kappa^+$ is regular then ${\rm Gal}(\mathscr{F},\kappa^+,\kappa^+)$ holds for every $\kappa$-complete filter $\mathscr{F}$ over $\kappa$.
\end{proposition}
\begin{proof}
Note that if $\kappa$ is boldface $\mathsf{GCH}$ and $\kappa^+$ is regular then every function $F\colon \kappa^+\rightarrow\mathcal{P}(\kappa)$ admits some $B\in\mathrm{Im}(F)$ such that $|(F^{-1})``B|=\kappa^+$.

Suppose $\langle A_\alpha\mid  \alpha < \kappa^+\rangle$ is a sequence of elements of a $\kappa$-complete filter $\mathscr{F}$ over $\kappa$. By the previous observation there is a $\kappa^+$-length subsequence which is constant. Hence the corresponding intersection  is this constant value, which belongs to $\mathscr{F}$. Thereby, ${\rm Gal}(\mathscr{F},\kappa^+,\kappa^+)$ holds.
\end{proof}
\begin{remark}
     Even ditching the regularity assumption on  $\kappa^+$ one can  show that ${\rm Gal}(\mathscr{F},\kappa,\kappa^+)$ holds for every $\kappa$-complete ultrafilter $\mathscr{F}$ over $\kappa.$ This is because every $F\colon \kappa^+\rightarrow\mathcal{P}(\kappa)$ admits a $B\in\mathrm{Im}(F)$ such that $|(F^{-1})``B|\geq \kappa.$
\end{remark} 

Anticipating the results in the next section we consider a natural generalization of Proposition~\ref{thmad}.
The crucial point is that both the family of sets with large intersection and the corresponding large set contained in all of them are explicitly constructed.
We will apply the claim to many pairs simultaneously, and since these objects are explicitly definable we do not need the axiom of choice in order to pick them.
\begin{definition}
We say that $\mathrm{Gal}_{\mathrm{Def}}(\mathscr{U},\lambda,\lambda)$ holds if for every sequence $\mathcal{C}=\l C_\alpha\mid\alpha<\lambda\r$, there is a subsequence $\l C_{\rho_j}\mid j<\lambda\r$  definable from $\mathcal{C}$ such that $\bigcap_{j<\lambda}C_{\rho_j}\in\mathscr{U}$.
\end{definition}

\begin{proposition} 
\label{clmgeneral} Suppose that $\kappa$ is measurable, $\kappa<\lambda$ is regular and there is no injection $f\colon \lambda\rightarrow \mathcal{P}(\kappa)$. Then, ${\rm Gal}_{\mathrm{Def}}(\mathscr{U},\lambda,\lambda)$ holds  
for every $\kappa$-complete ultrafilter $\mathscr{U}$ over $\kappa$.
\end{proposition}



If one seeks for a parallel of the above in \textsf{ZFC} then one may consider real-valued measurable cardinals, which can be described as measurables without the cardinal arithmetic.
A cardinal $\kappa$ is real-valued measurable if there exists a non-trivial $\kappa$-additive measure over $\kappa$.
The size of such cardinals, if exist at all, is at most $2^{\aleph_0}$.
Solovay proved that if there is a measurable cardinal $\kappa$ and one forces a $\kappa$-product of random reals then one obtains $2^{\aleph_0}=\kappa$ in the generic extension and $\kappa$ is real-valued measurable.
It is tempting to try an amalgamation between this theorem and the forcing construction of Abraham-Shelah \cite{MR830084}.
The framework will be similar, but rather than Cohen reals (added in the Abraham-Shelah model) one can try random reals.
The most difficult part is to verify that the \textit{Main lemma} \cite[Lemma 1.7]{MR830084} still holds true when replacing the Cohen's by the Random reals. This is in principle not evident for the original argument of \cite{MR830084} relied upon some specific properties of Cohen reals. 
Once this is accomplished, the failure of Galvin's property at $\aleph_1$ follows. However, there is the 
additional caveat of ensuring 
that the Abraham-Shelah poset 
preserves the real-valued measurability of $\kappa$. All of this suggests the following question:
\begin{question}
\label{qrealvalued} Is it consistent that $\kappa$ is a real-valued measurable cardinal and ${\rm Gal}(\mathscr{D}_{\aleph_1},\aleph_1,\kappa)$ fails?
\end{question}

\smallskip

\section{An application to ordinary partition relations}

Ramsey's theorem from \cite{MR1576401} says that $\omega\rightarrow(\omega,\omega)^2$.
That is, for every $c:[\omega]^2\rightarrow 2$ there exists an infinite monochromatic subset $A\subseteq\omega$.
A natural generalization is obtained by replacing $\omega$ with some cardinal $\kappa>\aleph_0$.
The resulting relation is  $\kappa\rightarrow(\kappa,\kappa)^2$ and implies that $\kappa$ is weakly compact.

There is yet another possible way to generalize Ramsey's theorem to uncountable cardinals and in this way one obtains a positive relation at small cardinals as well.
Recall that $\lambda\rightarrow(\kappa,\theta)^2$ means that for every $c:[\lambda]^2\rightarrow 2$, either there is $A\subseteq\lambda,|A|=\kappa$ such that $c``[A^2]=\{0\}$ or $B\subseteq\lambda, |B|=\theta$ such that $c``[B]^2=\{1\}$.
Ramsey's theorem generalizes to the statement $\kappa\rightarrow(\kappa,\omega)^2$ in which we increase the first component only.
A theorem of Erd\H{o}s, Dushnik and Miller says, indeed, that this relation holds at every infinite cardinal $\kappa$, see \cite{MR4862} and \cite{MR795592}.
In terms of graph theory this means that if $G$ is a complete graph of size $\kappa$ then either $G$ contains an edge-free subset of size $\kappa$ or an infinite clique.

Can one improve the positive relation $\lambda\rightarrow(\lambda,\omega)^2$?
The lightest possibility would be $\lambda\rightarrow(\lambda,\omega+1)^2$, but this relation does not hold at eve\-ry infinite cardinal $\lambda$ anymore.
Suppose that $\lambda>\cf(\lambda)=\omega$ and let $\lambda=\bigcup_{n<\omega}\Delta_n$ where $m\neq{n}\Rightarrow\Delta_m\cap\Delta_n=\emptyset$ and $|\Delta_n|<\lambda$ for every $n<\omega$.
Define $c:[\lambda]^2\rightarrow 2$ by letting $c(\alpha,\beta)=0$ iff there exists $n<\omega$ for which $\{\alpha,\beta\}\subseteq\Delta_n$.
A $0$-monochromatic set $A$ must be contained in some $\Delta_n$, so there is no such a set of size $\lambda$.
A $1$-monochromatic set $B$ satisfies $|B\cap\Delta_n|\leq 1$ for every $n<\omega$, so there is no such a set of order-type $\omega+1$.
Hence there is a class of cardinals which fail to satisfy $\lambda\rightarrow(\lambda,\omega+1)^2$.
On the other hand, if $\lambda=\cf(\lambda)>\aleph_0$ then $\lambda\rightarrow(\lambda,\omega+1)^2$, see \cite[Theorem 11.3]{MR795592}.
In fact, one can prove a slightly stronger statement.

\begin{proposition}
\label{thmstat} Suppose that $\lambda=\cf(\lambda)>\aleph_0$.
For every $c:[\lambda]^2\rightarrow 2$, either there is a stationary set $T\subseteq\lambda$ such that $c\upharpoonright[T]^2$ is $0$-monochromatic or there is $B\subseteq\lambda, {\rm otp}(B)=\omega+1$ such that $c\upharpoonright[B]^2$ is $1$-monochromatic.
\end{proposition}
\begin{proof}
Suppose that $c\colon [\lambda]^2\rightarrow 2$.
If there exists $B\subseteq\lambda$ of order type $\omega+1$ such that $c``[B]^2=\{1\}$ then we are done.
Suppose that there is no such $B$, and let $S$ be the set of limit ordinals of $\lambda$.
For every $\delta\in{S}$ choose a sequence $\bar{\alpha}_\delta = \langle\alpha^\delta_0,\ldots,\alpha^\delta_{n-1}\rangle$ such that $\bar{\alpha}_\delta^\frown\langle\delta\rangle$ is $1$-monochromatic and if $\max(\bar{\alpha}_\delta)<\xi<\delta$ then $\bar{\alpha}_\delta^\frown\langle\xi,\delta\rangle$ is not $1$-monochromatic.
The choice is possible by our assumption that there is no $1$-monochromatic sequence of length $\omega+1$.

By shrinking $S$ if needed, we may assume that $\ell{g}(\bar{\alpha}_\delta)=n$ for some fixed $n<\omega$ and every $\delta\in{S}$.
We remark that in this shrinking process we retain the fact that $S$ is stationary.
Let $\xi_\delta$ be the top-element of $\bar{\alpha}_\delta$ for every $\delta\in{S}$.
The function $h(\delta)=\xi_\delta$ is regressive on $S$, so by Fodor's lemma there is a stationary $T_0\subseteq{S}$ and a fixed ordinal $\xi<\lambda$ such that $\delta\in{T_0}\Rightarrow\xi_\delta=\xi$.
By repeating this process $n$-many times we obtain a stationary set $T_n$ and a fixed sequence $\bar{\alpha}$ such that $\bar{\alpha}^\frown\langle\delta\rangle$ is $1$-monochromatic and $\bar{\alpha}^\frown\langle\zeta,\delta\rangle$ is not $1$-monochromatic whenever $\delta\in{T_n}$ and $\max(\bar{\alpha})<\zeta<\delta$.

In particular, if $T=T_n\setminus (\max(\bar{\alpha})+1)$ then $T$ is a stationary subset of $\lambda$ and if $\zeta,\delta\in{T}, \zeta<\delta$ then necessarily $c(\zeta,\delta)=0$.
Otherwise, $\bar{\alpha}^\frown\langle\zeta\rangle$ will contradict the conclusion of the previous paragraph.
Thus, $c``[T]^2=\{0\}$ and we are done.
\end{proof}
As mentioned above, the statement of the proposition gives a bit more than just $\lambda\rightarrow(\lambda,\omega+1)^2$ since it yields a \emph{stationary} $0$-monochromatic set.
Notice that the argument applies to singular cardinals of uncountable cofinality for which the concepts of club and stationary subsets are sound.
However, if $T$ is a stationary subset of a singular cardinal $\lambda$ then it is possible that $|T|<\lambda$.
Therefore, one may wonder about $\lambda\rightarrow(\lambda,\omega+1)^2$ in such cardinals.
The following is \cite[Question 11.4]{MR795592}:

\begin{question}
\label{qehmr} Does the relation $\lambda\rightarrow(\lambda,\omega+1)^2$ hold for $\lambda>\cf(\lambda)>\omega$?
\end{question}

Actually, the question is phrased in \cite{MR795592} with respect to $\lambda=\aleph_{\omega_1}$, the first relevant instance.
We indicate that in \cite{MR795592} there appears a partial answer, based on canonization theorems of Shelah, which apply to strong limit singular cardinals.
Namely, if $\lambda>\cf(\lambda)>\omega$ and $\lambda$ is a strong limit cardinal then $\lambda\rightarrow(\lambda,\omega+1)^2$.
In some sense, this result gives many \textsf{ZFC} instances since for every $\kappa=\cf(\kappa)>\aleph_0$ there is a class of singular cardinals which are strong limit of cofinality $\kappa$.
On the other hand, for every specific $\lambda>\cf(\lambda)=\kappa$ one can choose $\theta<\lambda$ and force $2^\theta>\lambda$, thus locally it is not a theorem of \textsf{ZFC}.

A substantial progress with regard to the above question was made by Shelah in \cite{MR2494318}.
Using methods of pcf theory, Shelah proved that if $\lambda>\cf(\lambda)=\kappa>\aleph_0$ and $2^\kappa<\lambda$ then $\lambda\rightarrow(\lambda,\omega+1)^2$.
The assumption $2^\kappa<\lambda$ is a weakening of the assumption that $\lambda$ is a strong limit cardinal, but the overall question remains open if one wishes to eliminate any further assumption.

In this section we would like to replace $2^\kappa<\lambda$ by a version of Galvin's property. Our assumption on the Galvin property is also easily forced in the \textsf{ZFC} context if $\kappa=\cf(\lambda)$ is measurable.
This gives a slight improvement to Shelah's result, but it seems that the real importance of the $\textsf{ZFC}$ result is that it guides us towards the (tentative) direction of forcing the negative relation $\lambda\nrightarrow(\lambda,\omega+1)^2$. In effect, our approach indicates that one must kill all the pertinent instances of the Galvin property to obtain $\lambda\nrightarrow(\lambda,\omega+1)^2$.  As observed in the previous section, \textsf{AD} yields many instances of Galvin's property;  consequently, it also gives several instances of $\lambda\rightarrow(\lambda,\omega+1)^2$.

To this end, we must render the proof of \cite{MR2494318} by removing any use of choice apart from $\textsf{AC}_\omega$. Let us begin with models of \textsf{ZF}. Our first mission is to prove that $\kappa\rightarrow(\kappa,\omega+1)^2$ holds for many regular cardinals.
The proof of Theorem \ref{thmstat} seems to make use of the axiom of choice in two places.
Firstly, when one chooses the maximal green sequence $\bar{\alpha}_\delta$ below $\delta$ for every $\delta\in{S}$.
Secondly, when one employs Fodor's lemma (finitely many times).
The first issue is not a real obstacle, since finite sequences are well ordered even without choice.
The second issue is more substantial, but can be settled if one works with normal filters. In this respect,  we note that the club filter might not be normal, as proved by A. Karagila \cite{Karagila}. Measurability does not suffice, either. In effect, in \cite{NoNormal} Bilinsky and Gitik  constructed a model of \textsf{ZF} with a measurable cardinal which does not carry any normal ultrafilter. 

\begin{proposition}[\textsf{ZF}]
\label{clmadmeasurable} 
Suppose that $\mathscr{F}$ is a normal filter over $\kappa$ whose dual ideal contains $[\kappa]^{<\kappa}$. Then, for every $c\colon [\kappa]^2\rightarrow 2$ one can find either $A\in \mathscr{F}^+$, such that $c``[A]^2=\{0\}$ or $B\s \kappa$ with $\mathrm{otp}(B)=\omega+1$ such that $c``[B]^2=\{1\}.$  In particular, $\kappa\rightarrow(\kappa,\omega+1)^2$.
\end{proposition}
\begin{proof}
Let $c\colon [\kappa]^2\rightarrow 2$ be a coloring.
We refer to the first color as red and to the second color as green.
Repeat the arguments in the proof of Proposition~\ref{thmstat} with the following changes.
 Set $S:=\{\alpha<\kappa\mid \alpha\text{ is a limit ordinal}\}$. Clearly, $S\in\mathscr{F}^+$. 
For every $\delta\in{S}$ let $\bar{\alpha}_\delta$ be the first $<_{\rm lex}$-finite sequence which is green with $\delta$ and maximal with respect to this property.
Now use the normality and apply Fodor's lemma to $\mathscr{F}$ finitely many times to obtain a fixed maximal green sequence $\bar{\alpha}$ with respect to some set $A\in\mathscr{F}^+$.
It follows that $c``[A]^2=\{0\}$. Also, by our assumption on $\mathscr{F}$, $|A|=\kappa$. 
\end{proof}
Now we come to the main result of this section.

\begin{theorem}[\textsf{ZF}] 
\label{thm881ad} 
Suppose that $\aleph_0<\kappa=\cf(\lambda)<\lambda$ and  that $\kappa$ carries a normal measure $\mathscr{U}$. In addition, suppose that $\l\lambda_i\mid i<\kappa\r$ is a cofinal sequence in  $\lambda$ admitting a family  $\l \mathscr{U}_i\mid i<\kappa\r$ such that:
\begin{itemize}
   \item $\mathscr{U}_i$ is a normal filter whose dual ideal contains $[\lambda_i]^{<\lambda_i}$;
    \item $\mathrm{Gal}_{\mathrm{Def}}(\mathscr{U},\lambda_i,\lambda_i)$ holds.
\end{itemize}
Then, $\lambda\rightarrow(\lambda,\omega+1)^2$.
\end{theorem}
\begin{proof}
We may assume that for every $i<\kappa$, $\lambda_i>\kappa$.
Given a function $f\colon \kappa\rightarrow{\rm Ord}$ we define a rank function $\daleth(f)$ as follows.
Set $\daleth(f)=\alpha$ iff for every $\beta<\alpha$ one has $\daleth(f)\neq\beta$ and $\daleth(g)=\beta$ for some $g\colon \kappa\rightarrow{\rm Ord}$ which satisfies $g<_{\mathscr{U}}f$ (i.e., $\{\alpha<\kappa\mid g(\alpha)<f(\alpha)\}\in\mathscr{U}$).

Let $c\colon [\lambda]^2\rightarrow 2$ be a coloring and suppose that there is no $1$-monochromatic subset of order type $\omega+1$.
Define $\Delta_0:=\lambda_0$ and $\Delta_{1+i}:=[\lambda_i,\lambda_{i+1})$ for every $i<\kappa$.
By our assumption, there is a full-sized $0$-monochromatic subset of $\Delta_{i}$ for every $i<\kappa$.
Moreover, this set is explicitly constructible by the proof of Proposition~\ref{clmadmeasurable} (Indeed, for every $i$ we take the $<_{lex}$-minimal sequence $\vec{\alpha}^*_i$ such that the set of $\delta\in \Delta_i$ such that $\bar{\alpha}_\delta=\vec{\alpha}^*_i$ is in $\mathscr{U}_i$).
Hence we may assume without loss of generality that $c``[\Delta_i]^2=\{0\}$ for every $i<\kappa$.

For each $\alpha<\kappa$ let $\eta(\alpha)$ be the unique $i<\kappa$ for which $\alpha\in\Delta_i$.
For every $0<i<\kappa$ let $Seq_i$ be the set $\{\langle\alpha_0,\ldots,\alpha_{n-1}\rangle\mid \eta(\alpha_0)<\cdots<\eta(\alpha_{n-1})<i\}$.
For $i<\kappa$ and $\zeta\in\Delta_i$ we define a tree $\mathscr{T}_\zeta$ as follows.
For every $\bar{\alpha}=\langle\alpha_0,\ldots,\alpha_{n-1}\rangle\in Seq_i$ and every $\zeta\in\Delta_i$ we let $\bar{\alpha}\in\mathscr{T}_\zeta$ iff $\{\alpha_0,\ldots,\alpha_{n-1},\zeta\}$ is $1$-monochromatic under $c$.

By our assumption each $\mathscr{T}_\zeta$ is well-founded with respect to the reversed order.
Therefore, one can define a rank function $rk_\zeta$ over $Seq_i$, for every $\zeta\in\Delta_i$, by the following procedure.
If $\bar{\alpha}\in Seq_i-\mathscr{T}_\zeta$ then let $rk_\zeta(\bar{\alpha})=-1$.
If $\bar{\alpha}\in\mathscr{T}_\zeta$ then $rk_\zeta(\bar{\alpha})=\xi$ iff for every $\varepsilon\in\xi$ one has $rk_\zeta(\bar{\alpha})\neq\varepsilon$ and there exists an ordinal $\beta$ for which $rk_\zeta(\bar{\alpha}^\frown\langle\beta\rangle)=\varepsilon$.

The idea of this rank function is to express the degree of maximality of finite sequences below $\zeta$.
For a maximal $1$-monochromatic sequence below $\zeta$, $rk_\zeta$ assumes the value zero.
If there is more room for adding ordinals above $\max(\bar{\alpha})$ and keeping the $1$-monochromaticity with $\zeta$, then the rank grows.
Notice that for every $i<\kappa$ and every $\zeta\in\Delta_i$ the range of $rk_\zeta$ is bounded in $\lambda_{i+1}$ since it is an initial segment of $\lambda_{i+1}$.

Let $\Delta_i^{end}\subseteq\Delta_i$ be an end-segment such that for every $\bar{\alpha}\in Seq_i$ and $\gamma<\lambda_i$, if there is $\zeta\in \Delta_i^{end}$ such that $rk_\zeta(\bar{\alpha})=\gamma$ then there are unboundedly many such $\zeta$'s in $\Delta^{end}_i$ below $\lambda_{i+1}$.
The set $\Delta^{end}_i$ has a concrete definition: this end-segment is obtained by omitting bounded subsets of $\Delta_i$, which amounts in our case to the intersection of their complements.
Since $\lambda_{i+1}$ is regular, we obtain an end-segment of $\lambda_{i+1}$.

We define a set of triples $K$ as follows.
A triple $(\bar{\alpha},Z,f)$ belongs to $K$ iff $Z\in\mathscr{U}, f\colon \kappa\rightarrow{\rm Ord}$ and for some $0<i<\kappa, \bar{\alpha}\in Seq_i, \min(Z)>i$ and for every $j\in{Z}$ there exists $\zeta\in\Delta_j^{end}$ such that $rk_\zeta(\bar{\alpha})=f(j)$.
It is easy to verify that $K\neq\varnothing$.

Let $(\bar{\alpha}^*,Z^*,f^*)$ be a triple in $K$ for which $\daleth(f^*)$ is minimal amongst the elements of $K$.
Without loss of generality, all the elements of $Z^*$ are limit ordinals.
For every $j\in Z^*$ we isolate a section $\Sigma_j\subseteq\Delta_j^{end}$ by defining $\Sigma_j=\{\zeta\in\Delta_j^{end}\mid  rk_\zeta(\bar{\alpha}^*)=f^*(j)\}$.
Notice that $|\Sigma_j|=\lambda_{j+1}$.
Fix $j\in Z^*$.
We would like to understand what happens between the $j$th level and upper levels mentioned in $Z^*$.
For every $\ell\in Z^*$ such that $j<\ell$ and every $\zeta\in\Sigma_j$, let $\Sigma_j^\ell(\zeta)=\{\eta\in\Sigma_\ell\mid c(\zeta,\eta)=1\}$.
Let $L_\zeta=\{\ell\in Z^*\mid j<\ell,|\Sigma_j^\ell(\zeta)|=\lambda_{\ell+1}\}$, the set of \emph{large} $1$-monochromatic levels above $j$.
Similarly, let $T_\zeta=Z^*-L_\zeta$, the set of \emph{tiny} $1$-monochromatic levels above $j$.

Our goal is to garner many ordinals $\zeta$ with big $T_\zeta$, since then we will be able to remove the ``$1$" edges (there will be only a few of them) and create a $0$-monochromatic union of size $\lambda$.
We claim, therefore, that $L_\zeta\notin\mathscr{U}$ (and parallely, $T_\zeta\in\mathscr{U}$) for every $j\in Z^*,\zeta\in\Sigma_j$.
For suppose not.
Fix $i\in Z^*,\beta\in\Sigma_i$ such that $L_\beta\in\mathscr{U}$.
Define $\bar{\alpha}'={\bar{\alpha}}^{*\frown}\langle\beta\rangle, Z'=L_\beta$ and for $j\in Z'$ let $f'(j)=\min\{rk_\eta(\bar{\alpha}')\mid \eta\in\Sigma_i^j(\beta)\}$ and $f'(j)=0$ otherwise.
Notice that $(\bar{\alpha}',Z',f')\in{K}$.
We claim that $\daleth(f')<\daleth(f^*)$.
Indeed, for every $j\in Z'=L_\beta$ one has $f'(j)=rk_\eta(\bar{\alpha}')$ for some $\eta\in\Sigma_i^j(\beta)$, so $f'(j)=rk_\eta({\bar{\alpha}}^{*\frown}\langle\beta\rangle)<rk_\eta(\bar{\alpha}^*)=f^*(j)$.
Thus, $\daleth(f')<\daleth(f^*)$, contradicting the minimality of $\daleth(f^*)$.

We conclude, therefore, that $T_\zeta\in\mathscr{U}$ for every $j\in Z^*, \zeta\in\Sigma_j$.
Fix now $j\in Z^*$.

Let $\mathcal{A}^j=\l T_{\zeta}\mid \zeta\in\Sigma_j\r\subseteq \mathscr{U}$, apply $\mathrm{Gal}_{\mathrm{Def}}(\mathscr{U},\lambda_{j+1},\lambda_{j
+1})$ and find $M_j\subseteq \Sigma_j$ of cardinality $\lambda_{j+1}$  such that $\cap_{\zeta\in M_j}T_{\zeta}=:Y_j\in\mathscr{U}$.
We emphasize that we do not need the axiom of choice in order to define these objects, because of definability.
We render this process at every $j\in Z^*$.
Finally, let $Y=\Delta\{Y_j\mid j\in{Z^*}\}\in\mathscr{U}$.

For every $j\in{Y}$ we wish to define $A_j\subseteq M_j$ such that $|A_j|=\lambda_{j+1}$.
The sets $A_j$ will be $0$-monochromatic, and our goal is to show that their union is $0$-monochromatic as well.
We define these sets by induction on $j\in{Y}$.
So fix $j\in{Y}$ and define $A_j=\{\xi\in M_j\mid \forall i\in Y\cap j, \forall\zeta\in A_i, \xi\notin\Sigma_i^j(\zeta)\}$.
We claim that $|A_j|=\lambda_{j+1}$.
To see this, observe that if $i\in Y\cap{j}$ then for every $\zeta\in{A_i}$ the set $\Sigma_i^j(\zeta)$ is bounded in $\lambda_{j+1}$, hence $\lambda_{j+1}-\Sigma_i^j(\zeta)$ contains a final segment of $\lambda_{j+1}$.
Thus, by regularity of $\lambda_{j+1}$, $\bigcap\{(\lambda_{j+1}-\Sigma_i^j(\zeta))\mid \zeta\in{A_i}\}$ contains a final segment. Notice that  $A_j$ is obtained by intersecting this set with $M_j$, hence $|A_j|=\lambda_{j+1}$.

Define $A=\bigcup_{j\in Y}A_j$, so $|A|=\lambda$.
By proving that $c``[A]^2=\{0\}$ we will be done.
Pick $\alpha,\beta\in{A}$ such that $\alpha<\beta$.
If there exists $j\in\omega_1$ such that $\alpha,\beta\in{A_j}$ then $c(\alpha,\beta)=0$ since $A_j\subseteq A^j_0\subseteq\Sigma_j$ and $c``[\Sigma_j]^2=\{0\}$.
If not, then there are $i<j<\omega_1$ such that $\alpha\in{A_i}, \beta\in{A_j}$, and $i,j\in{Y}$.
By definition, $\beta\notin\Sigma_i^j(\alpha)$ and hence $c(\alpha,\beta)=0$ and the proof is accomplished.
\end{proof}

The above proof also gives the following corollary in \textsf{ZFC}, where the club filter over regular cardinals is always normal.
\begin{corollary}[ZFC] Suppose that $\mathscr{U}$ is a normal measure over $\kappa$, $\lambda$ is singular with $\cf(\lambda)=\kappa$ and that $\mathrm{Gal}(\mathscr{U},\lambda_i,\lambda_i)$ holds for cofinaly many regular cardinals below $\lambda$. Then, $\lambda\rightarrow(\lambda,\omega+1)^2$
\end{corollary}
The situation described in the previous corollary is a consequence of a small ultrafilter number: Suppose that $\kappa$ is measurable, and $\mathscr{U}$ is a normal ultrafilter over $\kappa$ with a base of $\mathscr{U}$ of size $\kappa^+$.
One says, in this case, that $\mathfrak{u}^{\rm nor}_\kappa=\kappa^+$.
It is possible to force $\mathfrak{u}^{\rm nor}_\kappa=\kappa^+$ even if $2^\kappa$ is arbitrarily large, see e.g. \cite{MR1632081} and \cite{MR3201820}.

It is easy to verify that if $\mathscr{U}$ witnesses $\mathfrak{u}^{\rm nor}_\kappa=\kappa^+$ then ${\rm Gal}(\mathscr{U},\mu,\mu)$ holds whenever $\mu=\cf(\mu)>\kappa^+$, see e.g. \cite{bgp}.

\begin{corollary}
\label{cormeasurable} If $\kappa$ is measurable and $\mathfrak{u}^{\rm nor}_\kappa=\kappa^+$ then $\lambda\rightarrow(\lambda,\omega+1)^2$ whenever $\kappa=\cf(\lambda)<\lambda$.
\end{corollary}

The corollary shows that the positive relation $\lambda\rightarrow(\lambda,\omega+1)^2$ may hold even if $2^\kappa>\lambda$. 
This fact was known already to Shelah and Stanley \cite{MR1782118}.
It is shown there that if $\lambda$ is a strong limit singular of uncountable cofinality then many versions of Cohen forcing preserve the positive relation $\lambda\rightarrow(\lambda,\omega+1)^2$; in particular, adding many Cohen subsets of $\kappa$ in a way that $2^\kappa>\lambda$.
Our corollary generalizes these results, since any $\kappa$-complete forcing notion will preserve the statement $\mathfrak{u}^{\rm nor}_\kappa=\kappa^+$.

In the above statements we used Galvin's property in order to prove $\lambda\rightarrow(\lambda,\omega+1)^2$.
The connection between ordinary partition relations and the structure of normal filters seems to be helpful in the opposite direction as well.
That is, from the assumption that $\lambda\rightarrow(\lambda,\omega+1)^2$ one can learn something about the Galvin property. 

\begin{proposition}
\label{clmctble} Let $\kappa=\cf(\kappa)>\aleph_0$ and let $\mathscr{F}$ be a normal filter over $\kappa$.
Suppose that $\neg_{\mathrm{st}}{\rm Gal}(\mathscr{F},\kappa,\kappa^+)$ is witnessed by $\mathcal{C}=\{C_\alpha\mid \alpha<\kappa^+\}$.
Then one can find $a=\{\alpha_n\mid n<\omega\}\subseteq\kappa^+$ such that $\bigcap\{(\kappa\setminus C_{\alpha_n})\mid n<\omega\}\neq\emptyset$.
\end{proposition}
\begin{proof} 
Assume toward contradiction that $\mathcal{C}=\{C_\alpha\mid \alpha<\kappa^+\}$ witnesses the strong failure of ${\rm Gal}(\mathscr{F},\kappa,\kappa^+)$, yet the complements are not overlapping in the sense that every infinite collection of them has empty intersection.
Define a coloring $c:[\kappa^+]^2\rightarrow 2$ as follows.
For $\alpha<\beta<\kappa^+$ let $c(\alpha,\beta)=0$ iff $\beta\in{C_\alpha}$.

Notice that there is no $1$-monochromatic sequence of length $\omega+1$ under $c$.
For if $\langle \alpha_n\mid n\leq\omega\rangle$ is such a sequence then $\alpha_\omega\in\bigcap\{(\kappa\setminus C_{\alpha_n})\mid n<\omega\}$, in contrary to our assumption at the beginning of the proof.
Likewise, there is no $0$-monochromatic set $A\in[\kappa^+]^{\kappa+\kappa}$.
For if $A=A_0\cup{A_1}$ with $\sup(A_0)<\min(A_1)$ and ${\rm otp}(A_0)={\rm otp}(A_1)=\kappa$, then $A_1\subseteq\bigcap\{C_\alpha\mid \alpha\in A_0\}$.
This is impossible since $\{C_\alpha\mid \alpha\in {A_0}\}\subseteq\mathcal{C}$.
Thus, the coloring $c$ witnesses $\kappa^+\nrightarrow(\kappa^+,\omega+1)^2$, which is an absurd since $\kappa^+$ is regular and uncountable.
\end{proof}

A similar argument becomes more powerful at successors of large cardinals.
By \cite{bgs} one can force a $\kappa$-complete ultrafilter $\mathscr{U}$ over a measurable cardinal $\kappa$ such that ${\rm Gal}(\mathscr{U},\kappa,\kappa^+)$ fails.

\begin{proposition}
\label{clmbgsfh} Suppose that $\kappa$ is measurable and $\mathscr{U}$ is a $\kappa$-complete ultrafilter over $\kappa$ for which $\neg_{\mathrm{st}}{\rm Gal}(\mathscr{U},\kappa,\kappa^+)$ is witnessed by the sequence $\mathcal{C}=\langle C_\alpha\mid \alpha<\kappa^+\rangle$. 
Then there is 
$S\in[\kappa]^\kappa$ such that $\bigcap\{(\kappa\setminus C_\alpha)\mid \alpha\in{S}\}$ is non-empty.
\end{proposition}
\begin{proof}
We define $c:[\kappa^+]^2\rightarrow 2$ as before, by letting $c(\alpha,\beta)=0$ iff $\beta\in{C_\alpha}$.
Of course, the definition is needed at $\alpha<\beta<\kappa^+$ only, since the coloring is symmetric.
Assume toward contradiction that there is no $S\in[\kappa]^\kappa$ such that $\bigcap\{(\kappa-C_\alpha)\mid \alpha\in{S}\}\neq\emptyset$.
By the same argument as in the previous claim, $c$ witnesses the negative relation $\kappa^+\nrightarrow(\kappa+\kappa,\kappa+1)^2$.
However, this relation holds if $\kappa$ is measurable, as proved in \cite{MR1968607}, a contradiction.
\end{proof}

The above statements show that there 
is a limitation on forcing empty intersection over the complements of sets which witness the strong failure of the Galvin property.
One can try, however, to force empty intersection only over certain families.
The following is tailored to the goal of obtaining a negative relation at singular cardinals.

Given $\lambda>\cf(\lambda)=\kappa>\aleph_0$ and a 
cofinal sequence $\langle \lambda_i\mid i<\kappa\rangle$ in $\lambda$, 
let $\Delta_0:=[0,\lambda_0)$ and $\Delta_{1+i}:=[\lambda_i,\lambda_{i+1})$.

\begin{proposition}
\label{clmnegdirection} Suppose that:
\begin{enumerate}
\item [$(\aleph)$] $\lambda>\cf(\lambda)=\kappa>\aleph_0$ and $\mathscr{F}$ is a normal filter over $\kappa$.
\item [$(\beth)$] $\{C_\alpha\mid \alpha<\lambda\}\subseteq\mathscr{F}$ witnesses $\neg_{\mathrm{st}}{\rm Gal}(\mathscr{F},\kappa,\lambda)$.
\item [$(\gimel)$] $\bigcap\{(\kappa\setminus C_{\alpha_n})\mid n<\omega\}=\emptyset$ whenever $\alpha_n\in\Delta_{i_n}$ for every $n<\omega$, and $m\neq n\Rightarrow \alpha_m\neq \alpha_n$.
\end{enumerate}
Then $\lambda\nrightarrow(\lambda,\omega+1)^2$.
\end{proposition}
\begin{proof}
For $\alpha<\beta<\lambda$, if there exists $i\in\kappa$ such that $\alpha,\beta\in\Delta_i$ then let $c(\alpha,\beta)=0$.
If $\alpha\in\Delta_i,\beta\in\Delta_j$ and $i<j$, let $c(\alpha,\beta)=0$ iff $\beta\in C_\alpha$.
One can verify that $c:[\lambda]^2\rightarrow 2$ witnesses $\lambda\nrightarrow(\lambda,\omega+1)^2$, so we are done.
\end{proof}

On a different path, Theorem~\ref{thm881ad} can be applied also to models of \textsf{AD}, where the distance between regular and measurable is quite short. This distance is even shorter under the extra assumption that $V=L(\mathbb{R})$, where every regular cardinal below $\Theta$ is measurable, as proved in \cite{MR2768698}.

\begin{corollary}[$\mathsf{AD}+V=L(\mathbb{R})$]
\label{thm881ad} 
Suppose that $\aleph_0<\kappa=\cf(\lambda)<\lambda$.
If  $\lambda$ is a limit of regular cardinals then $\lambda\rightarrow(\lambda,\omega+1)^2$.
\end{corollary}
\begin{proof}

For transparency, assume that $\kappa=\omega_1$.
Let $\lambda=\bigcup_{i<\omega_1}\lambda_i$, where $\langle \lambda_i\mid i<\omega_1\rangle$ is increasing and continuous, $\omega_1<\lambda_0$ is measurable and $\lambda_{i+1}$ is measurable for every $i<\omega_1$.
Let $\mathscr{U}$ denote the club filter over $\omega_1$, which is a normal ultrafilter under \textsf{AD}.
For every $i<\omega_1$ let $\mathscr{U}_{i+1}$ be the filter generated by the unbounded $\omega$-closed subsets of $\lambda_{i+1}$. This is a normal ultrafilter over $\lambda_{i+1}$ under \textsf{AD}+$V=L(\mathbb{R})$, as shown in \cite[Theorem~8.27]{MR2768698}. Also, $\mathrm{Gal}_{\mathrm{Def}}(\mathscr{U},\lambda_{i+1},\lambda_{i+1})$ holds by virtue of 
Proposition~\ref{clmgeneral}. Thereby, the assumptions of Theorem~\ref{thm881ad} are met. 
\end{proof}
It turns out that if $\lambda=\aleph_{\omega_1}$ we can do it  better:
\begin{corollary}[$\textsf{AD}$]
The relation $\aleph_{\omega_1}\rightarrow(\aleph_{\omega_1},\omega+1)^2$ holds.
\end{corollary}
\begin{proof}
For each limit $\alpha<\omega_1$ denote by  $\delta^1_\alpha$ the supremum of the lengths of all $\mathbf{\Delta}^1_\alpha$ prewellorderings.  It follows that, under \textsf{AD}, the sequence $\langle \mathbf{\delta}^1_{\alpha+1}\mid \alpha\in\mathrm{acc}(\omega_1)\rangle$ satisfies the following property: for each $\alpha\in \mathrm{acc}(\omega_1)$ the filter generated by the $\omega$-closed subsets of $\delta^1_{\alpha+1}$ is a normal measure \cite[\S5.4]{Jackson}. Also $\sup_{\alpha<\omega_1}\delta^1_{\alpha+1}=\aleph_{\omega_1}$. Now we can argue as in the previous corollary that  $\aleph_{\omega_1}\rightarrow(\aleph_{\omega_1},\omega+1)^2$.
\end{proof}
We finalize the current section with a $\textsf{ZFC}$ question:
\begin{question}
Does $\lambda\not\rightarrow(\lambda,\omega+1)^2$ have any consistency strength?
\end{question}

\subsection*{Acknowledgements}
We are deeply grateful to W. Chan for his wonderful insights and valuable feedback about the $\textsf{AD}$-related results.

\bibliographystyle{amsplain}
\bibliography{refer}

\end{document}